\documentclass[a4paper,twoside,11pt,reqno]{amsart}
\usepackage{fullpage}
\usepackage[T1]{fontenc}
\usepackage[utf8]{inputenc}
\usepackage{upgreek}
\usepackage{hyperref}
\usepackage{graphicx}
\usepackage{comment}
\usepackage[mathscr]{euscript}
\usepackage{booktabs}
\usepackage{mathtools}

\hypersetup{ocgcolorlinks=true,allcolors=testc}
\hypersetup{
     colorlinks   = true,
     citecolor    = black
}
\hypersetup{linkcolor=black}
\hypersetup{urlcolor=black}

\newcommand{\eqdef}{\stackrel{\mathrm{def}}{=}}

\newtheorem{theorem}{Theorem}
\newtheorem{lemma}[theorem]{Lemma}

\newtheorem{proposition}[theorem]{Proposition}
\newtheorem{corollary}[theorem]{Corollary}

\newtheorem*{definition*}{Definition}

\theoremstyle{remark}
\newtheorem{remark}[theorem]{Remark}

\newcommand*\diff{\mathop{}\!\mathrm{d}}
\newcommand{\R}{\mathbb{R}} 

\newcommand{\N}{\mathbb{N}}

\newcommand{\e}{\varepsilon}

\usepackage{dutchcal}

\renewcommand{\alpha}{\upalpha}
\renewcommand{\beta}{\upbeta}
\renewcommand{\gamma}{\upgamma}
\renewcommand{\delta}{\updelta}
\renewcommand{\zeta}{\upzeta}
\renewcommand{\eta}{\upeta}
\renewcommand{\theta}{\uptheta}
\renewcommand{\kappa}{\upkappa}
\renewcommand{\lambda}{\uplambda}
\renewcommand{\mu}{\upmu}
\renewcommand{\nu}{\upnu}
\renewcommand{\xi}{\upxi}
\renewcommand{\pi}{\uppi}
\renewcommand{\rho}{\uprho}
\renewcommand{\sigma}{\upsigma}
\renewcommand{\tau}{\uptau}
\renewcommand{\phi}{\upphi}
\renewcommand{\chi}{\upchi}
\renewcommand{\psi}{\uppsi}
\renewcommand{\omega}{\upomega}

\newcommand{\BB}{\mathsf{B}}

\newcommand{\E}{\mathbb{E}}

\newcommand{\mb}{\mathbb}

\newcommand{\mr}{\mathrm}

\newcommand{\vol}{\mathrm{vol}}

\renewcommand{\P}{\mathbb{P}}

\title{The dimension-free threshold for maximal sections of $\ell_p$}

\author{Silouanos Brazitikos}
\address{(S.~B.) Department of Mathematics and Applied Mathematics, University of Crete, 70013 Heraklion, Crete, Greece.}
\email{silouanb@uoc.gr} 

\author{Alexandros Eskenazis}
\address{(A.~E.) CNRS, Institut de Math\'ematiques de Jussieu, Sorbonne Universit\'e, Paris, France.}
\email{alexandros.eskenazis@imj-prg.fr}

\thanks{A.~E.~was supported by the European Research Council (ERC) under the European Union’s Horizon Europe programme, through the Starting Grant DiAnQuGe (grant agreement No. 101221523).}

\begin{document}

\maketitle
\vspace{-3mm}

\begin{abstract}
We prove the existence of a sharp threshold $p_0 \approx 26.265\ldots$, such that for every $n\in\mathbb{N}$ and $p\geq p_0$, the volume maximizing section of the unit ball of $\ell_p^n$ is the section with Ball's hyperplane $(\frac{e_1+e_2}{\sqrt{2}})^\perp$, whereas for $p\in(2,p_0)$ the asymptotic maximizer as the ambient dimension $n\to\infty$ is the hyperplane orthogonal to the diagonal vector $\frac{e_1+\cdots+e_n}{\sqrt{n}}$.
\end{abstract}

\bigskip

{\footnotesize

}


\section{Introduction}
Fix $p\in[1,\infty]$ and $n\in\N$. Let
\begin{equation}
\BB_p^n = \big\{x\in\R^n: \ \|x\|_p \leq 1\big\},
\end{equation}
be the closed unit ball of the normed space $\ell_p^n = (\R^n, \|\cdot\|_p)$, where for $p\in[1,\infty)$,
\begin{equation}
\forall \ x=(x_1,\ldots,x_n),\qquad  \|x\|_p = \big( |x_1|^p+\ldots+|x_n|^p\big)^{1/p}
\end{equation}
and $\|x\|_\infty = \max_{i=1,\ldots,n} |x_i|$.  The study of geometric parameters of these convex sets has a rich history,  with several precise computations remaining elusive for decades.  The purpose of this paper is to contribute to the following  long-standing question on sections of $\BB_p^n$. 

\smallskip
\noindent \emph{Question.}~For which unit vectors $a$ in $\R^n$ is the volume of $\BB_p^n\cap a^\perp$ maximal or minimal?
\smallskip

 This problem has its origins in the geometry of numbers and has served as a playground for a rich interplay of Fourier analytic and probabilistic techniques (see the survey \cite{NT22} and the references therein).  It was first Hadwiger \cite{Had72} and Hensley \cite{Hen79} who independently proved that the volume minimizing central sections of the unit cube $\BB_\infty^n$ are those with the hyperplanes $e_i^\perp$, which are parallel to the faces of the cube.  Conversely, the volume maximizing sections of $\BB_\infty^n$ were identified in the landmark paper \cite{Bal86} of Ball, who proved that for every unit vector $a$,
\begin{equation} \label{eq:ball}
\mr{vol}\big( \BB_\infty^n \cap a^\perp\big) \leq \mr{vol}\big( \BB_\infty^n \cap \big(\tfrac{e_1+e_2}{\sqrt{2}}\big)^\perp\big).
\end{equation}
We also refer to \cite{TT26} for a brilliant recent proof of Ball's cube slicing theorem.

The study of volume extremizing sections of $\BB_p^n$ for $p<\infty$ was initiated by Meyer and Pajor \cite{MP88},  who proved that the sections of $\BB_p^n$ with the coordinate hyperplanes $e_i^\perp$ have minimal volume for $p\geq2$ and maximal volume for $p\in[1,2]$. In the reverse direction, they showed that when $p=1$, the section of the cross-polytope $\BB_1^n$ with the hyperplane orthogonal to the main diagonal $\tfrac{e_1+\cdots+e_n}{\sqrt{n}}$ has minimal volume.  This was later extended to all values of $p\in[1,2]$ by Koldobsky \cite{Kol98} (see also \cite{ENT18} for a different probabilistic proof and \cite{Esk21,Mal25} for extensions).

The identification of volume maximizing sections of $\BB_p^n$ in the range $p\in(2,\infty)$ is substantially more subtle and has been raised a number of times \cite{Kol01, BN02, Ole03, KZ03, Kol05, MH11, Esk21,Tko21, NT22,Kon24,Kon25,GK26}.  Oleszkiewicz \cite{Ole03} proved that there exists a unique number $p_0 \approx 26.265\ldots$ with the following property: for every $2<p<p_0$, there exists $n=n(p)\in\N$ large enough for which the section of $\BB_p^n$ with the hyperplane $\big( \tfrac{e_1+\cdots+e_n}{\sqrt{n}}\big)^\perp$ has larger volume than the section with $\big(\tfrac{e_1+e_2}{\sqrt{2}}\big)^\perp$. In particular, one cannot expect a Ball-type extremal uniformly, for every $p>2$. The threshold $p_0$ is exactly the value of $p$ for which the volume of the section of the $p$-ball with Ball's hyperplane is equal, asymptotically as $n\to\infty$ and after proper normalization, with the diagonal section. Concretely, it is the unique solution of the equation 
\begin{equation} \label{eq:p0}
2^{2/p} \Gamma (\tfrac{1}{p})^3 = \pi  p^2 \Gamma(\tfrac{3}{p})
\end{equation}
in the interval $(2,\infty)$.  Delicate estimates on the dimension $n(p)$ were more recently obtained by K\"onig in \cite{Kon24}.  In the same work, Oleszkiewicz speculated that Ball-type hyperplanes may maximize the volume of sections for sufficiently large values of $p$. This was eventually confirmed by Eskenazis, Nayar and Tkocz in \cite{ENT24}, who proved that for every $n\in\N$,  $p\geq 10^{15}$ and every unit vector $a$, we have
\begin{equation} \label{eq:ent}
\mr{vol}\big( \BB_p^n \cap a^\perp\big) \leq \mr{vol}\big( \BB_p^n \cap \big(\tfrac{e_1+e_2}{\sqrt{2}}\big)^\perp\big).
\end{equation}

Identifying the volume maximizing sections of $\BB_p^n$ in the remaining interval $p\in(2,10^{15})$ is a formidable challenge, in part due to the expected transition of extremals around the value $p_0$ postulated by Oleszkiewicz's result. 
Our main result yields the sharp \emph{dimension-free} threshold for the maximal section problem in this missing range.

\begin{theorem} \label{thm:main}
Fix $n\geq2$ and let $p_0\approx 26.265\ldots$ be the solution of equation \eqref{eq:p0}.

\smallskip
\noindent {\bf (A)} If $p\geq p_0$, then for every unit vector $a$, we have
\begin{equation} \label{eq:goal-large}
\vol\big( \BB_p^n\cap a^\perp\big) \leq \mr{vol}\big( \BB_p^n \cap \big(\tfrac{e_1+e_2}{\sqrt{2}}\big)^\perp\big).
\end{equation}
\smallskip
\noindent {\bf (B)} If $p\in (2,p_0)$, then for every unit vector $a$, we have
\begin{equation} \label{eq:goal-small}
\frac{\vol\big( \BB_p^n\cap a^\perp\big)}{\vol\big( \BB_p^{n-1}\big)} \leq \limsup_{m\to\infty} \frac{\vol\big( \BB_p^m\cap\big(\frac{e_1+\ldots+e_m}{\sqrt{m}}\big)^\perp\big)}{\vol\big( \BB_p^{m-1}\big)}.
\end{equation}
\end{theorem}

The first part of Theorem \ref{thm:main} settles the maximal section problem for every $p$ above the critical threshold $p_0$, while the second part identifies the extremizers below the threshold in an \emph{asymptotic} fashion: after proper normalization of the section function, diagonals exhibit the worst behavior as the dimension tends to infinity. The same conclusion had previously been established by K\"onig in \cite{Kon25} for values $p\geq20$ which are far from even integers under the additional assumption that $\|a\|_\infty \leq \frac{1}{\sqrt2}$. As shall be explained below, removing these additional restrictions on $p$ and the normal vector $a$ is a delicate matter which requires a number of new ideas.

Theorem \ref{thm:main} is reminiscent of Haagerup's celebrated phase transition \cite{Haa81} for the optimal constants in the classical Khintchine inequality \cite{Khi22} which reads as follows. Let $A_p\in(0,1]$ be the largest constant such that for every $n\in\N$ and every $a_1,\ldots,a_n\in\R$, we have
\begin{equation} \label{eq:khintchine}
A_p \Big(\sum_{i=1}^n a_i^2\Big)^{1/2} \leq \bigg(\mb{E} \Big|\sum_{i=1}^n a_i \e_i\Big|^p \bigg)^{1/p},
\end{equation}
where $\e_1,\ldots,\e_n$ are independent Rademacher random variables. Then, there exists a threshold $p^* \approx 1.847\ldots$ such that $A_p = 2^{1/2-1/p}$ for $p\in[1,p^*]$ and $A_p= (\mb{E}|g|^p)^{1/p}$ for $p\in[p^*,2]$,  where $g$ is a standard Gaussian random variable. Equality holds when $a=\frac{e_1+e_2}{\sqrt{2}}$ for $p\in[1, p^*]$ and, asymptotically, when $a=\frac{e_1+\cdots+e_n}{\sqrt{n}}$ with $n\to\infty$ for $p\in[p^*,2]$. It is worth emphasizing that even in the study of the classical Khintchine inequality, the sharp constant in \eqref{eq:khintchine} for a \emph{fixed} number $n$ of Rademacher variables and $p\in(p^*,2)$ remains unknown. We refer to \cite{Saw85,LO95,OP00,KK01,BC02,BGMN05,Kon14,ENT18b,CKT21,CGT23,CST24,Kon25,RTW25,MRTT26} for other instances of such sharp transitions in the theory of optimal Khintchine-type inequalities and for further applications.


\section{Outline of the proof}


\subsection{Preliminaries} Fix $p\in[1,\infty)$, $n\in\N$ and consider the normalized section function of $\BB_p^n$, which we shall denote by $A_{n,p}:\mb{S}^{n-1}\to(0,\infty)$ and is given by
\begin{equation}
\forall \ a\in\mb{S}^{n-1},\qquad A_{n,p}(a) \eqdef \frac{\vol\big( \BB_p^n\cap a^\perp\big)}{\vol\big( \BB_p^{n-1}\big)}.
\end{equation}
In view of the symmetries of $\BB_p^n$ we shall always tacitly assume that the components of the normal vector $a\in\mb{S}^{n-1}$ satisfy $a_1\geq \cdots \geq a_n \geq0$. We start by presenting various known representations of the normalized section function, along with bounds which follow from them.

Consider the probability density $f_p:\R\to\R$ with $f_p(x) = e^{-\beta_p^p |x|^p}$, where \mbox{$\beta_p = 2\Gamma(1+\frac1p)$}. It is a classical fact due to \cite{MP88} (see also, e.g., \cite[Corollary~15]{ChasapisNayarTkocz2022}) that if $X_1,\ldots,X_n$ are i.i.d.~random variables with density $f_p$ and $g_a$ is the density of the weighted sum $\sum_{i=1}^n a_i X_i$, then the normalized section function can be expressed as
\begin{equation} \label{eq:dens-at-0}
\forall \ a\in\mb{S}^{n-1},\qquad A_{n,p}(a) = g_a(0).
\end{equation}
In other words, if we denote by $f_{p,a}(x) = \frac{1}{|a|} f_p(x/a)$ the density of $a X$, then
\begin{equation} \label{eq:conv}
\forall \ a\in\mb{S}^{n-1},\qquad A_{n,p}(a) = (f_{p,a_1} \ast \cdots \ast f_{p,a_n})(0).
\end{equation}
Since the densities $f_{p,a_i}$ are even, we can further express this as
\begin{equation} \label{eq:break-exp}
A_{n,p}(a) = \int_{\R} f_{p,a_1}(t) (f_{p,a_2} \ast \cdots \ast f_{p,a_n})(t)\, \diff t = \frac{1}{|a_1|} \mb{E} f_{p} \Big( \frac{1}{a_1} \sum_{i=2}^n a_i X_i\Big).
\end{equation}
This expression allows us to obtain bounds for the section function via H\"older's inequality.

\begin{lemma} \label{lem:holder}
Fix $p\in[1,\infty)$, $n\in\N$, $a\in\mb{S}^{n-1}$ with $1>a_1\geq\cdots\geq a_n\geq0$ and $X_1,\ldots,X_n$ independent random variables with density $f_p$. Denote by $\sigma = \sqrt{1-a_1^2}$. Then, if $g$ is the density of $B\eqdef\frac{1}{\sigma} \sum_{i=2}^n a_i X_i$ and $q\geq2$ and $q'\in[1,2]$ satisfy $\frac{1}{q}+\frac{1}{q'}=1$, we have
\begin{equation} \label{eq:holder}
A_{n,p}(a) \leq a_1^{-(1-1/q)} \sigma^{-1/q} q^{-1/pq} \|g\|_{2}^{2/q}.
\end{equation}
If $B'$ is an independent copy of $B$ and $\theta$ is the density of $\frac{B+B'}{\sqrt2}$, this is equivalent to
\begin{equation} \label{eq:holder-2}
A_{n,p}(a) \leq 2^{-1/2q} a_1^{-(1-1/q)} \sigma^{-1/q} q^{-1/pq} \theta(0)^{1/q}.
\end{equation}
In particular, choosing $q=\infty$, we get $A_{n,p}(a) \leq \frac{1}{a_1}.$
\end{lemma}

\begin{proof}
Renormalizing \eqref{eq:break-exp}, we can write
\begin{equation}
A_{n,p}(a) = \frac{1}{a_1} \mb{E} f_{p} \Big( \frac{\sigma}{a_1} \sum_{i=2}^n \frac{a_i}{\sigma} X_i\Big) = \frac{1}{a_1} \int_\R f_p\big(\tfrac{\sigma t}{a_1}\big) g(t) \, \diff t \leq \frac{1}{a_1} \big\| f_p\big( \tfrac{\sigma \ \cdot}{a_1}\big) \big\|_q \|g\|_{q'}.
\end{equation}
This yields the first bound since
\begin{equation*}
\big\| f_p\big( \tfrac{\sigma \ \cdot}{a_1}\big) \big\|_q^q = \frac{a_1}{\sigma} \int_\R e^{-q\beta_p^p |x|^p} \,\diff x = \frac{a_1 q^{-1/p}}{\sigma}
\end{equation*}
and $\|g\|_{q'} \leq \|g\|_2^{2/q}$ because $g$ is a probability density. Finally, observe that, since $g$ is even,
\begin{equation} \label{eq:symm-426}
\|g\|_2^2 = \sqrt{2} \int_\R g\big(\sqrt{2} x\big)^2 \, \diff x = \tfrac{1}{\sqrt{2}} \ (\widetilde{g}\ast \widetilde{g})(0) = \tfrac{1}{\sqrt2}\theta(0),
\end{equation}
where $\widetilde{g}(x)=\sqrt{2} g(\sqrt{2} x)$ is the density of $B/\sqrt2$, and this proves \eqref{eq:holder-2}.
\end{proof}

In \cite[Section~3.3]{Bar01}, Barthe used a similar expression for the section function along with the Brascamp--Lieb inequality \cite{BL76,Bal89} to obtain the following bound, which will be useful later on.

\begin{proposition} [Barthe] \label{prop:barthe}
Fix $p\in[2,\infty)$ and $n\in\N$. Then, for every $a\in\mb{S}^{n-1}$, we have
\begin{equation} \label{eq:barthe}
A_{n,p}(a) \leq \bigg( \prod_{i=1}^n (1-a_i^2)^{-(1-a_i^2)} \bigg)^{1/2-1/p}.
\end{equation}
\end{proposition}

A different representation of the section function, going back to \cite{MP88} (see also \cite[Theorem~3.2]{Kol05} or \cite{Kon25}), can be obtained via the Fourier transform. Let $\gamma_p:\R\to\R$ be the function
\begin{equation} \label{eq:gammap}
\forall \ s\in\R,\qquad \gamma_p(s) \eqdef \frac{1}{\Gamma(1+\frac{1}{p})} \int_0^\infty \cos(sx) e^{-x^p} \, \diff x.
\end{equation}
An application of the Fourier inversion formula on \eqref{eq:conv} gives
\begin{equation} \label{eq:section-fourier}
\forall \ a\in\mb{S}^{n-1}, \qquad A_{n,p}(a) = \tfrac{2}{\pi}\Gamma\big(1+\tfrac1p\big) \int_0^\infty \prod_{i=1}^n \gamma_p(a_i s)\,\diff s.
\end{equation}
When $p=\infty$, we have $\gamma_\infty(s) = \frac{\sin s}{s}$. In analogy with Ball's proof of the cube slicing theorem \cite{Bal86} (see also \cite{ENT24b} for extensions), consider the function $\Psi_p:[1,\infty)\to\R_+$ given by
\begin{equation} \label{eq:Psip}
\forall \ u\geq1, \qquad \Psi_p(u) \eqdef \sqrt{u} \int_0^\infty \big| \gamma_p(s)\big|^u\, \diff s.
\end{equation}
A similar function appeared in Haagerup's proof \cite{Haa81} of the optimal constant in the Khintchine inequality. We derive the following bound for $A_{n,p}$ in terms of the Ball--Haagerup function.

\begin{lemma} \label{lem:holder-Psi}
Fix $p\in[1,\infty)$ and $n\in\N$. Then, for every $a\in\mb{S}^{n-1}$, we have
\begin{equation} \label{eq:holder-Psi}
A_{n,p}(a) \leq  \tfrac{2}{\pi}\Gamma\big(1+\tfrac1p\big) \prod_{i=1}^n \Psi_p\big(a_i^{-2}\big)^{a_i^2}.
\end{equation}
Modulo permutations and changes of sign, equality holds when $a=\frac{ e_1+\cdots+ e_k}{\sqrt{k}}$ for some even $k\leq n$.
\end{lemma}

\begin{proof}
By H\"older's inequality, since $\sum_{i=1}^n a_i^2=1$,
\begin{equation}
\tfrac{\pi}{2} \Gamma\big(1+\tfrac{1}{p}\big)^{-1} A_{n,p}(a) \stackrel{\eqref{eq:section-fourier}}{\leq} \prod_{i=1}^n \bigg( \int_0^\infty \big| \gamma_p(a_is)\big|^{a_i^{-2}} \, \diff s \bigg)^{a_i^2} = \prod_{i=1}^n \bigg( |a_i|^{-1} \int_0^\infty \big| \gamma_p(\sigma)\big|^{a_i^{-2}} \, \diff \sigma \bigg)^{a_i^2},
\end{equation}
which is the desired inequality. The equality cases also follow readily.
\end{proof}

Bounds on the Ball--Haagerup function $\Psi_p$ play a crucial role in the proof of Theorem \ref{thm:main}. For the moment, we record a simple estimate which will be useful later on.

\begin{lemma} \label{lem:HY}
Fix $p\in[1,\infty)$. Then, for every $u\geq2$, we have
\begin{equation} \label{eq:HY}
\Psi_p(u) \leq \frac{\pi}{2\Gamma(1+\frac{1}{p})} \bigg( \frac{u}{u-1}\bigg)^{(u-1)(\frac12-\frac1p)}.
\end{equation}
\end{lemma}

\begin{proof}
First, observe that since $f_p$ is even,
\begin{equation}
\gamma_p(s) = \frac{1}{2\Gamma(1+1/p)} \int_\R e^{-|x|^p} e^{isx} \, \diff x = \widehat{f_p}\big( \tfrac{\beta_p s}{2\pi}\big).
\end{equation}
Therefore, Beckner's Hausdorff--Young inequality with an optimal constant \cite{Bec75} yields
\begin{equation}
\Psi_p(u) = \tfrac12\sqrt{u} \|\gamma_p\|_u^u = \frac{\pi}{\beta_p} \sqrt{u} \|\widehat{f_p}\|_u^u \leq \frac{\pi}{2\Gamma(1+\frac1p)} \sqrt{u} \bigg( \frac{u'^{1/u'}}{u^{1/u}} \bigg)^{u/2} \|f_p\|_{u'}^u,
\end{equation}
where $\frac{1}{u}+\frac{1}{u'}=1$, which is the required bound since $\|f_p\|_{u'} = u'^{-1/pu'}$.
\end{proof}

We shall also use a probabilistic representation of the normalized section function via negative moments. Variants of this formula for general convex bodies can be traced to the work \cite{KK05} of Kalton and Koldobsky. In the case of $\BB_p^n$, Chasapis, Nayar and Tkocz \cite[Lemma~14]{ChasapisNayarTkocz2022} proved that if $X_1,\ldots,X_n$ are i.i.d.~random variables with density $f_p$, then
\begin{equation} \label{eq:negative-moments}
\forall \ a\in\mb{S}^{n-1},\qquad A_{n,p}(a) = \lim_{s\downarrow -1} \frac{1+s}{2} \mb{E} \Big| \sum_{i=1}^n a_i X_i\Big|^s.
\end{equation}
Using \eqref{eq:negative-moments}, it was shown in \cite[Proposition 6]{ENT24} that if $R_1,\ldots,R_n$ are i.i.d.~random variables with density $\alpha_p^{-1} x^p e^{-x^p} {\bf 1}_{x>0}$, where $\alpha_p = \frac{1}{p}\Gamma(1+\frac{1}{p})$, and $\xi_1,\ldots,\xi_n$ are i.i.d.~random vectors uniformly distributed on the unit sphere $\mb{S}^2$, independent of the random variables $R_i$, then
\begin{equation} \label{eq:spherical}
\forall \ a\in\mb{S}^{n-1},\qquad A_{n,p}(a) = \Gamma\big(1+\tfrac{1}{p}\big) \mb{E} \Big| \sum_{i=1}^n a_i R_i\xi_i \Big|^{-1},
\end{equation}
where $|\cdot|$ stands for the Euclidean norm $\|\cdot\|_2$. This spherical negative moment identity
can also be derived directly from
\eqref{eq:section-fourier} through the Fourier representation of
negative moments of Gorin and Favorov
\cite[Lemma~3]{GorinFavorov1990}. Indeed, for a random vector
\(Z\) in \(\mathbb R^3\), we have
\begin{equation} \label{eq:g-f}
\E|Z|^{-1}=\frac1{2\pi^2}\int_{\mathbb R^3}\E e^{i\langle t,Z\rangle} \, \frac{\diff t}{|t|^2},
\end{equation}
whenever the integral is well-defined. Moreover, if $\xi$ is uniformly distributed on $\mb{S}^2$, then the random variable $\langle v, \xi\rangle$ is uniformly distributed on $[-|v|,|v|]$ for every $v\in\R^3$ by Archimedes' hat-box theorem. Therefore, for $u\in\R^3$, 
\begin{equation} \label{eq:F-tr}
\begin{split}
\E_{\xi,R} &e^{i\langle u,R\xi\rangle}
= \mb{E} \Big[ \frac{1}{2R|u|} \int_{-R|u|}^{R|u|} e^{i\theta} \, \diff \theta \Big] =
\E \ \frac{\sin(R|u|)}{R|u|} \\ &  = \frac{p}{\Gamma(1+\frac1p)} \int_0^\infty \frac{\sin(r|u|)}{|u|} r^{p-1} e^{-r^p}\,\diff r 
 = \frac{1}{\Gamma(1+\frac1p)} \int_0^\infty \cos(r|u|) e^{-r^p} \, \diff r \stackrel{\eqref{eq:gammap}}{=}
\gamma_p(|u|),
\end{split}
\end{equation}
where the penultimate equality follows from integration by parts. Using independence and integration in polar coordinates, we thus get
\begin{equation*}
\Gamma\big(1+\tfrac{1}{p}\big) \mb{E} \Big| \sum_{i=1}^n a_i R_i\xi_i\Big|^{-1} \stackrel{\eqref{eq:g-f}}{=} \frac{\Gamma(1+\frac1p)}{2\pi^2} \int_{\R^3} \prod_{j=1}^n \mb{E} e^{i \langle a_j t, R\xi\rangle} \, \frac{\diff t}{|t|^2} \stackrel{\eqref{eq:F-tr}}{=} \tfrac{2}{\pi} \Gamma\big(1+\tfrac1p\big) \int_0^\infty \prod_{i=1}^n \gamma_p(a_i s) \, \diff s
\end{equation*}
and this indeed yields \eqref{eq:spherical} because of \eqref{eq:section-fourier}.

Finally, another classical application of Archimedes' hat-box theorem (see, e.g., \cite[Lemma~24]{ChasapisNayarTkocz2022}) shows that if $X, Y$ are two independent rotationally invariant random vectors in $\R^3$, then
\begin{equation} \label{eq:min}
\mb{E} |X+Y|^{-1} = \mb{E} \min\big\{|X|^{-1},|Y|^{-1}\big\}.
\end{equation}
Using this  identity, one can derive the following estimate for the section function.
\begin{lemma} \label{lem:min-bound}
Fix $p\in[1,\infty)$ and $n\in\N$. If $R_1,\ldots,R_n$ are i.i.d.~random variables with density $\alpha_p^{-1} x^p e^{-x^p} {\bf 1}_{x>0}$, where $\alpha_p = \frac{1}{p}\Gamma(1+\frac{1}{p})$, and $\xi_1,\ldots,\xi_n$ are i.i.d.~random vectors uniformly distributed on the unit sphere $\mb{S}^2$, independent of the random variables $R_i$, then for every unit vector $a$, we have
\begin{equation} \label{eq:min-bound}
A_{n,p}(a) \leq \Gamma\big(1+\tfrac1p\big) \ \mb{E}\Big[ \min \Big\{\big|a_1 R_1 \xi_1 + a_2 R_2 \xi_2 \big|^{-1}, \mb{E} \Big| \sum_{i\geq3} a_i R_i \xi_i \Big|^{-1} \Big\} \Big].
\end{equation}
\end{lemma}

\begin{proof}
Consider the independent rotationally invariant random vectors $X= a_1 R_1 \xi_1 + a_2 R_2 \xi_2$ and $Y=\sum_{i\geq3} a_i R_i \xi_i $. Then, we have
\begin{equation}
\Gamma\big(1+\tfrac1p\big)^{-1} A_{n,p}(a) \stackrel{\eqref{eq:spherical}}{=}  \mb{E}|X+Y|^{-1} \stackrel{\eqref{eq:min}}{=} \mb{E}_{X,Y} \min \big\{|X|^{-1},|Y|^{-1}\big\} \leq \mb{E} \min\big\{|X|^{-1},\mb{E}|Y|^{-1}\big\},
\end{equation}
by Jensen's inequality for the concave function $t\mapsto \min\{t,a\}$, where $a\in\R$.
\end{proof}

This bound was crucially used in \cite{ChasapisNayarTkocz2022} and \cite{ENT24} to prove the stability and resilience respectively of Ball's cube slicing theorem \cite{Bal86}. It shall play a central role in our proof as well.

To conclude this discussion, we will rewrite Theorem \ref{thm:main} in terms of the normalized section function. For $p\geq p_0$, the desired inequality \eqref{eq:goal-large} becomes
\begin{equation} \label{eq:goal-large-2}
A_{n,p}(a) \leq A_{n,p}\big(\tfrac{e_1+e_2}{\sqrt{2}}\big) \stackrel{\eqref{eq:negative-moments}}{=} A_{2,p}\big(\tfrac{e_1+e_2}{\sqrt{2}}\big) = \frac{1}{\|(\frac{1}{\sqrt{2}},\frac{1}{\sqrt{2}})\|_p} = 2^{\frac12-\frac1p}.
\end{equation}
For $p\in(2,p_0)$, the central limit theorem allows us to rewrite inequality \eqref{eq:goal-small} as
\begin{equation} \label{eq:goal-small-2}
\begin{split}
A_{n,p}(a) & \leq \lim_{m\to\infty} A_{m,p}\Big(\frac{e_1+\cdots+e_m}{\sqrt{m}}\Big) \stackrel{\eqref{eq:spherical}}{=} \Gamma\big(1+\tfrac1p\big) \lim_{m\to\infty} \mb{E}\Big| \frac{1}{\sqrt{m}} \sum_{i=1}^m R_i \xi_i\Big|^{-1}
\\ & = \Gamma\big(1+\tfrac1p\big) \big(\mb{E} R^2 / 3\big)^{-1/2} \mb{E} |G|^{-1} = \bigg( \frac{6}{\pi} \frac{\Gamma(1+\frac1p)^3}{\Gamma(1+\frac3p)} \bigg)^{1/2},
\end{split}
\end{equation}
where $G$ is a standard Gaussian in $\R^3$. Observe that the two upper bounds of \eqref{eq:goal-large-2} and \eqref{eq:goal-small-2} match exactly at $p=p_0$, so Theorem \ref{thm:main} can be equivalently reformulated as
\begin{equation} \label{eq:final-goal}
A_{n,p}(a) \leq \max\bigg\{ 2^{\frac12-\frac1p}, \bigg( \frac{6}{\pi} \frac{\Gamma(1+\frac1p)^3}{\Gamma(1+\frac3p)} \bigg)^{1/2} \bigg\}.
\end{equation}
Observe that, as shown in \cite{Ole03}, the maximum coincides with the first term if and only if $p\geq p_0$.


\subsection{Outline of the proof} The proof of Theorem \ref{thm:main} draws ideas from the aforementioned works of Eskenazis, Nayar, Tkocz \cite{ENT24} and K\"onig \cite{Kon25}, combined with several new insights. A key ingredient is a sharp upper bound for the Ball--Haagerup special function $\Psi_p$ of \eqref{eq:Psip}.

\begin{theorem} \label{thm:Psi}
Fix $p\geq2$. Then, for every $u\geq2$, we have
\begin{equation} \label{eq:Psi-bound}
\Psi_p(u) \leq \max\big\{ \Psi_p(2), \Psi_p(\infty)\big\}.
\end{equation}
\end{theorem}

Theorem \ref{thm:Psi} was previously proven by K\"onig in \cite{Kon25} for values $p\geq20$ which are far from even integers, while the case $p=\infty$ coincides with Ball's integral inequality \cite{Bal86}. K\"onig's proof laboriously implements the ingenious method of Nazarov and Podkorytov \cite{NP00}, who used distribution functions to establish similar inequalities for the special functions of \cite{Haa81} and \cite{Bal86}. In contrast, our proof of Theorem \ref{thm:Psi} crucially relies on a delicate polynomial envelope for the Fourier transform $\gamma_p$ (see Section \ref{sec:Qp}). This is more reminiscent of the hands-on approach towards the integral inequalities found in the original papers \cite{Haa81,Bal86}. Our polynomial bound is then combined with the Hausdorff--Young inequality and a monotonicity argument to yield the conclusion \eqref{eq:Psi-bound}.

Combining the sharp bound \eqref{eq:Psi-bound} for the special function $\Psi_p$ with Lemma \ref{lem:holder-Psi}, we deduce that if a unit vector $a$ satisfies $\|a\|_\infty\leq \frac{1}{\sqrt{2}}$, then each $a_i^{-2}\geq2$ and thus
\begin{equation} \label{eq:derivation-small-coefficients}
\begin{split}
A_{n,p}(a) \leq  \tfrac{2}{\pi}\Gamma\big(1+\tfrac1p\big) \prod_{i=1}^n \Psi_p\big(a_i^{-2}\big)^{a_i^2} & \stackrel{\eqref{eq:Psi-bound}}{\leq} \tfrac{2}{\pi}\Gamma\big(1+\tfrac1p\big) \max\big\{\Psi_p(2),\Psi_p(\infty)\big\} 
\\ & = \max\bigg\{ 2^{\frac12-\frac1p}, \bigg( \frac{6}{\pi} \frac{\Gamma(1+\frac1p)^3}{\Gamma(1+\frac3p)} \bigg)^{1/2} \bigg\},
\end{split}
\end{equation}
where the last identity follows from the equality cases of Lemma \ref{lem:holder-Psi} and the computations in \eqref{eq:goal-large-2}, \eqref{eq:goal-small-2}. This is the desired conclusion \eqref{eq:final-goal} of Theorem \ref{thm:main}. To circumvent the constraint on $\|a\|_\infty$, several additional ideas are required. The strategy is different depending on the value of $p$, so we describe the main ideas separately. Assume that $\|a\|_\infty > \frac{1}{\sqrt2}$.

\smallskip

\noindent {\bf (1)} In the range $2<p<4$, the proof uses only the Brascamp--Lieb bound \eqref{eq:barthe}.

\smallskip

\noindent {\bf (2)} In the range $4\leq p \leq 26$, we start by applying the H\"older estimate \eqref{eq:holder} and symmetrize as in \eqref{eq:holder-2}. The resulting bound involves the density at 0 of the weighted sum $(B+B')/\sqrt2$ of two independent copies of $B=\frac{1}{\sigma} \sum_{i=2}^n a_i X_i$. Since the associated coefficient vector has $\ell_\infty$-norm at most $\frac{1}{\sqrt{2}}$, the Fourier estimate \eqref{eq:derivation-small-coefficients} then applies and gives \eqref{eq:final-goal}.

\smallskip

\noindent {\bf (3)} In the range $p>26$, we start by considering the same random vector $B$ but we need to further distinguish various cases depending on the distribution of its coefficient vector. 

\noindent $\bullet$ If all components of the vector are sufficiently small, we use the probabilistic representation \eqref{eq:break-exp} and replace the random vector $B$ by an appropriate Gaussian, while carefully tracking the error of approximation and the size of the corresponding Gaussian integral.

\noindent $\bullet$ If some components have intermediate size, then we apply the H\"older estimate \eqref{eq:holder} and symmetrize via \eqref{eq:holder-2} as in the range $4\leq p\leq 26$. However, a mere application of the Fourier estimate \eqref{eq:derivation-small-coefficients} is not sufficient here. Instead, we shall need a stronger form of Theorem \ref{thm:Psi} which includes a \emph{deficit term} and shall be proven in Section \ref{sec:Psi} below (see Theorem \ref{thm:Psi-deficit}). 

\noindent $\bullet$ Finally, if the coefficient vector is sufficiently spiked, we resort to the negative moment representation \eqref{eq:spherical} and to the inductive method of \cite{ENT24}. While the main ideas in this range are similar to those of \cite{ENT24}, more delicate technical work is needed to allow for $p$ to be as small as 26 instead of the previously known bound of $10^{15}$.

\smallskip

To recapitulate, our proof of Theorem \ref{thm:main} relies on two conceptually distinct refinements of previous methods \cite{ENT24,Kon25}. First, we circumvent distribution function arguments and obtain a direct proof of the sharp upper bound \eqref{eq:Psi-bound} for the Ball--Haagerup function, including a deficit term, thus obtaining the conclusion \eqref{eq:final-goal} when $\|a\|_\infty\leq \frac{1}{\sqrt{2}}$. In the range $\|a\|_\infty > \frac{1}{\sqrt{2}}$, two morally different treatments are needed. When $a$ is far from the extremizer, we estimate the section function using probabilistic representations \eqref{eq:holder} and the Brascamp--Lieb inequalities \eqref{eq:barthe}, often after a symmetrization argument, to derive the conclusion \eqref{eq:final-goal}. These serve as more efficient substitutes to the stability of Ball's cube slicing theorem \cite{MR22,ChasapisNayarTkocz2022} that was used in \cite{ENT24}. Finally, when $a$ is near the extremizer, a careful implementation of the inductive argument of Eskenazis, Nayar and Tkocz \cite{ENT24} completes the proof.


\section{A polynomial envelope for $\gamma_p$} \label{sec:Qp}

This section contains some technical estimates which we shall need later on for the Fourier transform $\gamma_p$. The interested reader may wish to skip the proofs of this section during the first reading of the paper, as they consist of mostly elementary yet tedious manipulations.

Fix $p\in[2,\infty)$ and recall that the function $\gamma_p$ is given by
\begin{equation} \label{eq:gammap-again}
\forall \ s\in\R,\qquad \gamma_p(s) = \frac{1}{\Gamma(1+\frac{1}{p})} \int_0^\infty \cos(sx) e^{-x^p} \, \diff x.
\end{equation}
Expanding the cosine for small $s$, we have
\begin{equation*}
\gamma_p(s) = \frac{1}{\Gamma(1+\frac{1}{p})} \int_0^\infty  e^{-x^p} \, \diff x - \frac{s^2}{2 \Gamma(1+\frac{1}{p})} \int_0^\infty x^2 e^{-x^p} \, \diff x + O(s^4) = 1 - \frac{\Gamma(1+\frac3p)}{6\Gamma(1+\frac1p)} \ s^2+ O(s^4)
\end{equation*}
as $s\to0$. Denote the absolute value of the coefficient of $s^2$ by $c_p$, namely
\begin{equation} \label{eq:cp}
c_p \eqdef \frac{\Gamma(1+\frac3p)}{6\Gamma(1+\frac1p)}.
\end{equation}
K\"onig proved in \cite[Proposition~3.1]{Kon25} that for $p\geq5$ and $s\in[0,3]$,
\begin{equation} \label{eq:Konig-envelope}
\gamma_p(s) \leq e^{-c_p s^2}.
\end{equation}
Observe that the Taylor expansions of the two functions coincide up to second order. We shall need a polynomial improvement of \eqref{eq:Konig-envelope} which captures the size of $\gamma_p$ up to fourth order. To this end, observe that we have
\begin{equation*}
\gamma_p(s) = 1 - c_p s^2 + \frac{s^4}{24 \Gamma(1+\frac{1}{p})} \int_0^\infty x^4 e^{-x^p} \, \diff x + O(s^6) =  1 - c_p s^2 + \frac{\Gamma(1+\frac5p)}{120 \Gamma(1+\frac{1}{p})} \ s^4 + O(s^6) 
\end{equation*}
as $s\to0$. Consider the parameter $\theta_p$ given by
\begin{equation} \label{eq:thetap}
\theta_p \eqdef 1 - \frac{3}{5} \frac{\Gamma(1+\frac1p)\Gamma(1+\frac5p)}{\Gamma(1+\frac3p)^2}
\end{equation}
and, for $p>2$, the function $Q_p:[0, (\theta_pc_p)^{-1/2}]\to [0,1]$ given by
\begin{equation}
\forall \ s\in\big[0, (\theta_pc_p)^{-1/2}\big],\qquad Q_p(s) \eqdef \big( 1-\theta_p c_p s^2\big)^{1/\theta_p}.
\end{equation}
Observe that the Taylor expansions of $\gamma_p$ and $Q_p$ around $s=0$ match up to fourth order and, trivially, $Q_p(s) \leq e^{-c_p s^2}$ for every $s$.

We start by proving the following simple properties of $p\mapsto\theta_p$.

\begin{lemma} \label{lem:thetap}
The map $p\mapsto\theta_p$ is increasing for $p>2$ and $a\mapsto\theta_{1/a}$ is concave on $(0,\frac12)$.
\end{lemma}

\begin{proof}
Let $\sigma(a) \eqdef \frac{\Gamma(1+a)\Gamma(1+5a)}{\Gamma(1+3a)^2}$ and we will prove that $\sigma$ is increasing and convex on $(0,\frac12)$. First, by the series representation \cite[6.3.16]{AS64} for the derivative of the polygamma function,
\begin{equation}
(\log\sigma)'(a) = \sum_{k=1}^\infty \bigg( \frac{6}{k+3a} - \frac{5}{k+5a}-\frac{1}{k+a}\bigg) \geq 0
\end{equation}
by the concavity of the function $s\mapsto \frac{s}{k+sa}$, which means that $\log\sigma$ is increasing and so is $\sigma$. For the second claim, we shall prove that $\log\sigma$ is convex, which readily implies that $\sigma$ is also convex. Using the integral formula \cite[6.4.1]{AS64} for the derivatives of the polygamma function,
\begin{equation}
(\log\sigma)''(a) = \ \int_0^\infty \bigg( \frac{1}{e^{v}-1} + \frac{1}{e^{v/5}-1} - \frac{2}{e^{v/3}-1} \bigg) v e^{-av} \, \diff v
\end{equation}
and the conclusion follows by the convexity of the function $x\mapsto \frac{1}{e^{1/x}-1}$, where $x>0$.
\end{proof}

The main technical result of this section is an estimate refining K\"onig's bound \eqref{eq:Konig-envelope}.

\begin{proposition} \label{prop:Qp}
Fix $p>2$. Then, for every $s\in[0, c_p^{-1/2}]$, we have
\begin{equation} \label{eq:Qp-bound}
\gamma_p(s) \leq Q_p(s).
\end{equation}
Moreover, if $p\geq 3$, inequality \eqref{eq:Qp-bound} holds for every $s\in[0,3]$.
\end{proposition}

\begin{proof}
As $\gamma_p(0)=Q_p(0)=1$, it suffices to show that $\gamma_p/Q_p$ is nonincreasing. This is equivalent to the nonpositivity of the function
\begin{equation}
\zeta_p(s) \eqdef (1-\theta_p c_p s^2) \gamma_p'(s) + 2 c_p s \gamma_p(s)
\end{equation}
for $s$ in the corresponding interval. Expanding the cosine function, we can write
\begin{equation}
\gamma_p(s) = \frac{1}{\Gamma(1+\frac{1}{p})} \int_0^\infty \cos(sx) e^{-x^p} \, \diff x = \sum_{n=0}^\infty \frac{(-1)^n \Gamma(\frac{2n+1}{p}) s^{2n} }{(2n)! \Gamma(\frac{1}{p})}
\end{equation}
If $a_{n} \eqdef \frac{\Gamma(\frac{2n+1}{p})}{(2n)!\Gamma(\frac1p)}$ are the alternating coefficients in the expansion above, we can expand $\zeta_p$ as
\begin{equation} \label{eq:taylor}
    \begin{split}
         \zeta_p(s) &  = (1-\theta_p c_p s^2) \sum_{n=1}^\infty (-1)^n\cdot 2n\cdot  a_n s^{2n-1} + 2c_p \sum_{n=0}^\infty (-1)^n a_n s^{2n+1}
        \\ &= - 2\sum_{n=0}^\infty (-1)^n (n+1) a_{n+1} s^{2n+1} - 2\theta_p c_p \sum_{n=0}^\infty (-1)^n n a_n s^{2n+1} + 2c_p \sum_{n=0 }^\infty (-1)^n a_n s^{2n+1}
        \\ & = 2\sum_{n=0}^\infty (-1)^n \big(c_p a_n - \theta_p c_p n a_n - (n+1) a_{n+1}\big)  s^{2n+1}.
    \end{split}
\end{equation}
{\bf Claim.} For every $n\in\N$,
\begin{equation}
\kappa_n \eqdef c_p a_n - \theta_p c_p n a_n - (n+1) a_{n+1} \leq 0.
\end{equation}

\begin{proof}
Substituting the values of $a_n, c_p = a_1, \theta_p = 1-\frac{2a_2}{a_1^2}$, we want to prove
\begin{equation} \label{eq:an-ineq}
\kappa_n=(1-n)a_1a_n + \frac{2na_2 a_n}{a_1} - (n+1) a_{n+1} \leq0.
\end{equation}
Now consider the sequence
\begin{equation*}
\rho_n \eqdef \frac{n a_n}{a_{n-1}} = \frac{1}{2(2n-1)} \frac{\Gamma(\frac{2n+1}{p})}{\Gamma(\frac{2n-1}{p})} = \frac{1}{2p} \frac{\Gamma(\frac{2n+1}{p})}{\Gamma(\frac{2n-1}{p}+1)} = \frac{1}{2p\Gamma(1-\frac2p)} \int_0^1 t^{\frac{2n+1}{p}-1} (1-t)^{-\frac2p} \, \diff t,
\end{equation*}
where the last equality is the relation of the $B$ and $\Gamma$ functions. Setting $u=t^{2/p}$, we can further express this as
\begin{equation} \label{eq:moment}
\rho_n = \int_0^1 u^{n-1} \, \diff\mu(u),
\end{equation}
where $\diff\mu(u) = \frac{u^{1/2} ( 1-u^{p/2})^{-2/p}}{4\Gamma(1-\frac2p)}\, \diff u$ is a positive measure. Dividing by $a_n$, we can rewrite \eqref{eq:an-ineq} as
\begin{equation}
0\leq \rho_{n+1}- n\rho_2  + (n-1) \rho_1 = \int_0^1 t^{n} -1  -n( t-1) \, \diff \mu(t)
\end{equation}
which is true since the integrand is pointwise nonnegative.
\end{proof}

Grouping the terms of the expansion \eqref{eq:taylor} in pairs, we can rewrite this as
\begin{equation} \label{eq:zeta}
\frac{ \zeta_p(s)}{2s}  = \sum_{n=0}^\infty (-1)^n \kappa_n s^{2n} = \sum_{n=0}^\infty \big( \kappa_{2n} - \kappa_{2n+1} s^2\big) s^{4n}. 
\end{equation}
We will prove that each $\kappa_{2n}-\kappa_{2n+1}s^2$ is nonpositive. By the claim above, it suffices to check this at the right endpoint of the interval. Expressing everything in terms of $\rho_n$,
\begin{equation*}
\begin{split}
\kappa_{2n}&-\kappa_{2n+1} c_p^{-1} = - a_{2n} \big( \rho_{2n+1} - \rho_1 - 2n(\rho_2-\rho_1) \big) + \frac{a_{2n+1}}{c_p} \big( \rho_{2n+2} - \rho_1 - (2n+1)(\rho_2-\rho_1) \big)
\\ & = - a_{2n} \Big( \rho_{2n+1} - \rho_1 - 2n(\rho_2-\rho_1) \big) - \frac{\rho_{2n+1}}{(2n+1)c_p} \big( \rho_{2n+2} - \rho_1 - (2n+1)(\rho_2-\rho_1) \big) \Big) 
\\ & \stackrel{\eqref{eq:moment}}{=} - a_{2n} \int_0^1 u^{2n} - 1  - 2n (u-1) - \frac{\rho_{2n+1}}{(2n+1)c_p} \big( u^{2n+1} - 1 - (2n+1)(u-1) \big) \, \diff \mu(u).
\end{split}
\end{equation*}
Consider the polynomial $p_m(u) = u^m - 1 - m(u-1)$ which is nonnegative on $[0,1]$. We shall next show that the ratio $p_{m+1}/p_m$ is increasing in $(0,1)$. Indeed, we have
\begin{equation}
\begin{split}
p_{m+1}'(u) p_m(u) & - p_{m+1}(u) p_m'(u) = u^{2m} - m^2 u^{m+1} + 2(m^2-1) u^m -m^2 u^{m-1}+1 
\\ & = \big(1-u^m\big)^2 - m^2 u^{m-1} (1-u)^2
\end{split}
\end{equation}
and this is nonnegative if and only if
\begin{equation}
\Big( \frac{1-u^m}{1-u}\Big)^2 = \bigg( \sum_{k=0}^{m-1} u^k \bigg)^2 = \sum_{r=0}^{2m-2} \sigma_r u^r \geq m^2 u^{m-1},
\end{equation}
where $\sigma_r = \min\{r+1,2m-1-r\}$ is the number of pairs $(i,j)\in\{0,1,\ldots,m-1\}^2$ with $i+j=r$. Using the AM-GM inequality $u^r + u^{2m-2-r} \geq 2 u^{m-1}$ and $\sigma_r = \sigma_{2m-2-r}$, we get
\begin{equation}
\sum_{r=0}^{2m-2} \sigma_r u^r \geq \sum_{r=0}^{2m-2} \sigma_r  u^{m-1} = m^2 u^{m-1}
\end{equation}
and this proves the monotonicity of $p_{m+1}/p_m$. Hence, for $n\geq1$,
\begin{equation}
\frac{p_{2n+1}(u)}{p_{2n}(u)} \leq \lim_{v\to 1} \frac{p_{2n+1}(v)}{p_{2n}(v)} = \frac{p_{2n+1}''(1)}{p_{2n}''(1)} = \frac{2n+1}{2n-1}.
\end{equation}
This gives
\begin{equation}
\begin{split}
\kappa_{2n}-\kappa_{2n+1} c_p^{-1} &= -a_{2n} \int_0^1 \Big( 1 - \frac{\rho_{2n+1}}{(2n+1)c_p} \frac{p_{2n+1}(u)}{p_{2n}(u)}\Big) p_{2n}(u) \, \diff\mu(u)
\\ & \leq -a_{2n} \Big( 1 - \frac{\rho_{2n+1}}{(2n-1)c_p} \Big) \int_0^1 p_{2n}(u) \, \diff\mu(u).
\end{split}
\end{equation} 
To conclude, it suffices to observe that $\rho_{2n+1} \leq (2n-1)c_p$ since \eqref{eq:moment} gives $\rho_{2n+1} \leq \rho_1 = c_p$. The case $n=0$ is trivial since $\kappa_0=\kappa_1=0$.

For the second claim of the proposition, the same reasoning shows that it suffices to check that each coefficient of the expansion \eqref{eq:zeta} is nonpositive for $s=3$ and
\begin{equation}
\kappa_{2n} - 9\kappa_{2n+1} \leq -a_{2n} \Big( 1 - \frac{9\rho_{2n+1}}{2n-1} \Big) \int_0^1 p_{2n}(u) \, \diff\mu(u).
\end{equation}
If $n\geq2$, then
\begin{equation}
\frac{9\rho_{2n+1}}{2n-1} \leq 3\rho_1 =  3c_p = \frac{\Gamma(1+\frac3p)}{2\Gamma(1+\frac1p)} \leq \tfrac{1.4}{1.6} < 1,
\end{equation}
where we used that $\min_{x\in[1,2]} \Gamma(x) > 0.8$ and $\max_{x\in[1,5/2]}\Gamma(x) = \Gamma(5/2) < 1.4$. The case $n=0$ is again trivial since $\kappa_0=\kappa_1=0$. Finally, for $n=1$, we have to prove
\begin{equation} \label{eq:k2k3}
\kappa_2-9\kappa_3
\leq0.
\end{equation}
As $u^2-1-2(u-1)=(1-u)^2$, $u^3-1-3(u-1)=(1-u)^2(u+2)$ and $9a_3=3a_2\int u^2\, \diff\mu$,
\begin{equation}
\kappa_2-9\kappa_3=a_2\Big(3\int_0^1u^2\, \diff\mu\int_0^1(1-u)^2(u+2)\, \diff\mu-\int_0^1(1-u)^2\, \diff\mu\Big).
\end{equation}
Hence the desired inequality is equivalent to
\begin{equation}\label{eq:goal-gamma}
3\int_0^1u^2\, \diff\mu\cdot\Big(2+\frac{\int_0^1u(1-u)^2\, \diff\mu}{\int_0^1(1-u)^2\, \diff\mu}\Big)<1 .
\end{equation}
The density of $(1-u)^2\, \diff\mu$ is proportional to
\begin{equation}
\Big(\frac{1-u^{p/2}}{1-u}\Big)^{-2/p} \cdot u^{1/2}(1-u)^{2-\frac2p}.
\end{equation}
The first factor, which we shall denote by $\sigma(u)$, is decreasing in $u$ since $u\mapsto u^{p/2}$ is convex on $[0,1]$ and thus $u\mapsto\frac{1-u^{p/2}}{1-u}$ is increasing. Therefore, Chebyshev's rearrangement inequality for the law $\pi$ of the $\mathrm{Beta}(\frac32,3-\frac2p)$ distribution, yields
\begin{equation}
\frac{\int_0^1u(1-u)^2\, \diff\mu}{\int_0^1(1-u)^2\, \diff\mu} = \frac{\int_0^1 u \sigma(u)\,\diff\pi(u)}{\int_0^1 \sigma(u)\,\diff\pi (u)} \leq \int_0^1 u\,\diff\pi(u) = \frac{\frac32}{\frac32+3-\frac2p}=\frac{3}{9-\frac4p}\le\frac9{23}.
\end{equation}
On the other hand,
\begin{equation}
\int_0^1u^2\, \diff\mu=\frac{3a_3}{a_2}=\frac1{14}\cdot\frac{\Gamma(1+\frac7p)}{\Gamma(1+\frac5p)} .
\end{equation}
As a function of $a=\frac1p\in[0,\frac13]$, the ratio $\Gamma(1+\frac7p)/\Gamma(1+\frac5p)$ is log-convex, since
\begin{equation}
\frac{\diff^2}{\diff a^2}\log\frac{\Gamma(1+7a)}{\Gamma(1+5a)}=\sum_{k=1}^\infty \bigg(\frac{49}{(k+7a)^2}-\frac{25}{(k+5a)^2}\bigg)\ge0,
\end{equation}
each term being nonnegative because $7(k+5a)\ge5(k+7a)$. Therefore the ratio is maximal at an endpoint. It equals $1$ at $a=0$, and at $a=\frac13$, we have $\frac{\Gamma(\frac{10}3)}{\Gamma(\frac83)}<1.93.$
Hence $\int_0^1u^2\, \diff\mu<\frac{1.93}{14}$. This completes the proof of inequality \eqref{eq:goal-gamma} since $3\cdot\frac{1.93}{14}\cdot (2+\frac{9}{23}) < 1$.
\end{proof}

Finally, we record the following consequence of Proposition \ref{prop:Qp} which will be useful later on. For $p>2$, consider the parameter $\lambda_p$ given by
\begin{equation} \label{eq:lambdap}
\lambda_p\eqdef
\begin{cases}
	\frac12,&2<p<3,\\[1mm]
	\frac1{3\Gamma(1+\frac1p)},&p\geq3.
\end{cases}
\end{equation}

\begin{corollary} \label{cor:Qp}
Fix $p>2$. Then, for $s>0$,
\begin{equation}
|\gamma_p(s)|>\lambda_p \quad\Longrightarrow\quad 0<\gamma_p(s)\le Q_p(s).
\end{equation}
\end{corollary}

\begin{proof}
Assume first that $2<p<3$ and $|\gamma_p(s)|>\frac12$. Integrating by parts twice, we get 
\begin{equation} \label{eq:two-bp}
\begin{split}
\Gamma\big(1+\tfrac1p\big) |\gamma_p(s)| & = \frac{p}{s} \bigg| \int_0^\infty \sin(sx) x^{p-1} e^{-x^p}\,\diff x\bigg| = \frac{p}{s^2}\bigg| \int_0^\infty \cos(sx) \big(x^{p-1} e^{-x^p}\big)'\,\diff x\bigg| 
\\ & \leq \frac{p}{s^2} \int_0^\infty \big| \big(x^{p-1} e^{-x^p}\big)' \big| \, \diff x = \frac{2p}{s^2} \left(\frac{p-1}{pe}\right)^{\frac{p-1}{p}},
\end{split}
\end{equation}
where the last equality follows since the function $x\mapsto x^{p-1}e^{-x^p}$ is increasing on $[0,x^*]$ and decreasing on $[x^*,\infty)$, where $x^* = (1-\frac1p)^{1/p}$. 

We shall now prove the one-dimensional estimate
\begin{equation} \label{eq:1D}
\frac{2pc_p}{\Gamma(1+\frac1p)} \left(\frac{p-1}{pe}\right)^{\frac{p-1}{p}} \le\frac12.
\end{equation}
Indeed, let $a=\frac1p\in(\frac13,\frac12)$ and consider
\begin{equation}
h(a) \eqdef \log \frac{pc_p}{\Gamma(1+\frac1p)} = -\log a + \log\Gamma(1+3a) - 2\log\Gamma(1+a) -  \log 6.
\end{equation}
Using the classical formula \cite[6.4.10]{AS64} for the derivatives of the polygamma function,
\begin{equation}
h''(a) = \frac{1}{a^2} + \sum_{k=1}^\infty \bigg(\frac{9}{(k+3a)^2} - \frac{2}{(k+a)^2} \bigg) > \sum_{k=0}^\infty \frac{9}{(k+\frac32)^2} - \sum_{k=1}^\infty \frac{2}{(k+\frac13)^2} > 0
\end{equation}
on $(\frac13,\frac12)$ and thus $h$ is convex. Therefore, the graph of $h$ is below the chord joining the points $(\frac13,h(\frac13))$, $(\frac12,h(\frac12))$ and this gives
\begin{equation}
\frac{pc_p}{\Gamma(1+\frac1p)} \leq \Big(\frac{1}{\sqrt{\pi}}\Big)^{\frac6p-2}\Big(\frac{1}{2\Gamma(\frac43)^2}\Big)^{3-\frac6p}.
\end{equation}
Thus,
\begin{equation}
\frac{2pc_p}{\Gamma(1+\frac1p)} \left(\frac{p-1}{pe}\right)^{\frac{p-1}{p}} \leq 2\pi^{1-\frac3p}\,\big(2\,\Gamma(\tfrac43)^2\big)^{\frac6p-3}\Big(\frac{p-1}{pe}\Big)^{\frac{p-1}{p}}.
\end{equation}
The resulting upper bound is easily seen to be log-convex in $a=\frac1p$ so its maximum value is at one of the endpoints $\frac13, \frac12$, both of which are below $\frac12$. This concludes the proof of \eqref{eq:1D}, which combined with \eqref{eq:two-bp} and the assumption on $\gamma_p$ implies that $s<c_p^{-1/2}$. Therefore the domination of $\gamma_p$ by $Q_p$ follows from Proposition \ref{prop:Qp}. 

To prove that $\gamma_p$ is positive on this interval use the inequality $\cos x \geq 1-\frac{x^2}{2}$ to get
\begin{equation} \label{eq:lb-gammap}
\gamma_p(s) \geq 1- c_p s^2
\end{equation}
which is evidently nonnegative for $s<c_p^{-1/2}$.

Now let $p\ge3$. Integrating by parts once like in \eqref{eq:two-bp} gives
\begin{equation} \label{eq:bp-1-time}
\Gamma\big(1+\tfrac1p\big)|\gamma_p(s)| \le \frac{p}{s}\int_0^\infty x^{p-1} e^{-x^p}\,\diff x = \frac1{s}
\end{equation}
and thus $|\gamma_p(s)|>\lambda_p$ implies $s<3$. The domination of $\gamma_p$ by $Q_p$ follows from Proposition \ref{prop:Qp}.

It remains only to prove that $\gamma_p(s)\ge0$ for $s\in[0,3]$. This was previously proven by K\"onig in \cite[Lemma~4.2]{Kon25} for $p\geq15$. Using the inequality
\begin{equation}
\cos x \ge 1-\frac{x^2}{2} +\frac{x^4}{24} -\frac{x^6}{720}
\end{equation}
and integrating against $e^{-x^p}\diff x$, we get
\begin{equation}
\gamma_p(s) \ge 1-c_ps^2+ \frac{\Gamma(1+\frac5p)}{120\Gamma(1+\frac1p)}s^4 - \frac{\Gamma(1+\frac7p)}{5040\Gamma(1+\frac1p)}s^6.
\end{equation}
The polynomial $P(s)$ on the right is non-negative on $[0,3]$ for every $p\ge3$. To see this, first calculate
\begin{equation}
P'(s)=-\frac{s}{840\,\Gamma(1+\frac1p)}\Big(280\,\Gamma\big(1+\tfrac3p\big)-28\,\Gamma\big(1+\tfrac5p\big)s^2+\Gamma\big(1+\tfrac7p\big)s^4\Big).
\end{equation}
The expression in the parenthesis is a quadratic polynomial in $s^2$ with positive leading coefficient, whose discriminant is negative if and only if
\begin{equation}
\Gamma\big(1+\tfrac5p\big)^2<\tfrac{10}{7}\,\Gamma\big(1+\tfrac3p\big)\Gamma\big(1+\tfrac7p\big).
\end{equation}
This holds since $\log\Gamma$ is convex, thus $\Gamma(1+\frac5p)^2\leq\Gamma(1+\frac3p)\Gamma(1+\frac7p)$. Therefore $P'(s)<0$ for every $s>0$ and $P(s)\geq P(3)$ for $s\in[0,3]$. 

Finally, we shall prove that
\begin{equation} \label{eq:P(3)}
P(3) = 1-\frac{3\,\Gamma(1+\frac3p)}{2\,\Gamma(1+\frac1p)}+\frac{27\,\Gamma(1+\frac5p)}{40\,\Gamma(1+\frac1p)}-\frac{81\,\Gamma(1+\frac7p)}{560\,\Gamma(1+\frac1p)} \geq0.
\end{equation}
Write $a=\frac1p\in(0,\frac13]$ and $R_k(a)=\frac{\Gamma(1+ka)}{\Gamma(1+a)}$ for $k\geq1$, so that
\begin{equation} \label{eq:PR}
P(3)=1-\tfrac32R_3(a)+\tfrac{27}{40}R_5(a)-\tfrac{81}{560}R_7(a).
\end{equation}
Each $R_k$ is log-convex, and in particular convex, on $[0,\infty)$ since, by \cite[6.4.10]{AS64},
\begin{equation}
\frac{\diff^2}{\diff a^2}\log R_k(a)=k^2(\log\Gamma)''(1+ka)-(\log\Gamma)''(1+a)=\sum_{j=1}^\infty \bigg(\frac{k^2}{(j+ka)^2}-\frac{1}{(j+a)^2}\bigg)\geq0,
\end{equation}
as $k(j+a)\geq j+ka$. Partition $[0,\frac13]$ into the intervals $[0,\frac1{12}]$, $[\frac1{12},\frac16]$, $[\frac16,\frac14]$, $[\frac14,\frac7{24}]$ and $[\frac7{24},\frac13]$. On each of these intervals, denote by $\mathrm{chord}(R_k)$, $k\in\{3,7\}$, the affine function which coincides with $R_k$ at the endpoints and by $\mathrm{tang}(R_5)$ the tangent line of $R_5$ at $\frac1{24}$ on the first interval and at the right endpoint on the other four.
By convexity, $R_3\leq\mathrm{chord}(R_3)$, $R_7\leq\mathrm{chord}(R_7)$ and $R_5\geq\mathrm{tang}(R_5)$ on each interval, thus
\begin{equation} \label{eq:affine}
P(3)\geq1-\tfrac32\,\mathrm{chord}(R_3)+\tfrac{27}{40}\,\mathrm{tang}(R_5)-\tfrac{81}{560}\,\mathrm{chord}(R_7).
\end{equation}
The right-hand side of \eqref{eq:affine} is piecewise affine in $a$, hence it suffices to check that it is positive at the two endpoints of each interval, where the chords coincide with $R_3$ and $R_7$. The relevant values of $R_3$ and $R_7$ are
\begin{center}
\begin{tabular}{c|cccccc}
\toprule
$a$ & $0$ & $\frac1{12}$ & $\frac16$ & $\frac14$ & $\frac7{24}$ & $\frac13$\\
\midrule
$R_3(a)$ & $1$ & $0.94586$ & $0.95527$ & $1.01397$ & $1.06083$ & $1.11985$\\
$R_7(a)$ & $1$ & $0.93056$ & $\frac76$ & $1.77444$ & $2.31327$ & $\frac{28}9$\\
\bottomrule
\end{tabular}
\end{center}
and the values of $R_5$ and $R_5'$ at the tangency points are
\begin{center}
\begin{tabular}{c|ccccc}
\toprule
$a$ & $\frac1{24}$ & $\frac16$ & $\frac14$ & $\frac7{24}$ & $\frac13$\\
\midrule
$R_5(a)$ & $0.93699$ & $1.01394$ & $\frac54$ & $1.43697$ & $1.68489$\\
$R_5'(a)$ & $-0.82638$ & $1.90467$ & $3.86273$ & $5.16056$ & $6.80842$\\
\bottomrule
\end{tabular}
\end{center}
where the entries $R_7(\frac16)=\frac{\Gamma(13/6)}{\Gamma(7/6)}=\frac76$, $R_5(\frac14)=\frac{\Gamma(9/4)}{\Gamma(5/4)}=\frac54$ and $R_7(\frac13)=\frac{\Gamma(10/3)}{\Gamma(4/3)}=\frac{28}9$ are exact and $R_3(\frac13)=\Gamma(\frac43)^{-1}$. Plugging these values in, the right-hand side of \eqref{eq:affine} takes the following values at the endpoints of the five intervals
\begin{center}
\begin{tabular}{c|ccccc}
\toprule
interval & $[0,\frac1{12}]$ & $[\frac1{12},\frac16]$ & $[\frac16,\frac14]$ & $[\frac14,\frac7{24}]$ & $[\frac7{24},\frac13]$\\
\midrule
left endpoint & $0.0111$ & $0.0239$ & $0.0248$ & $0.0472$ & $0.0200$\\
right endpoint & $0.0558$ & $0.0828$ & $0.0661$ & $0.0441$ & $0.0075$\\
\bottomrule
\end{tabular}
\end{center}
all of which are positive. Therefore $P(3)>0$ for every $a\in(0,\frac13]$, that is, for every $p\geq3$.
\end{proof}


\section{Bounds for the Ball--Haagerup function} \label{sec:Psi}

In this section, we present our optimal bounds for the Ball--Haagerup function $\Psi_p$ of Theorem \ref{thm:Psi}, asserting that for every $p\geq2$  and $u\geq2$,
\begin{equation} \label{eq:Psi-bound-again}
\Psi_p(u) \leq \max\big\{ \Psi_p(2), \Psi_p(\infty)\big\}.
\end{equation}
The case $p=2$ is trivial as $\gamma_2(s)=e^{-s^2/4}$ and $\Psi_2(u)\equiv \sqrt{\pi}$ for every $u\geq1$.  For $p>2$, consider the parameter $v_p$ given by
\begin{equation} \label{eq:vp}
v_p\eqdef
\begin{cases}
	\frac{2p}{p-2},&2< p<3,\\[1mm]
	3,& 3\leq p\leq6,\\[1mm]
	2,&p>6.
\end{cases}
\end{equation}
To prove Theorem \ref{thm:Psi} for $p>2$, we shall first establish \eqref{eq:Psi-bound-again} for small $u\in[2,v_p]$ using the Hausdorff--Young bound \eqref{eq:HY}. Then, we will use the domination of $\gamma_p$ by $Q_p$ of Proposition \ref{prop:Qp} to identify a monotone quantity which is larger than $\Psi_p(u)/\sqrt{u}$ and \eqref{eq:Psi-bound-again} will follow.

\begin{lemma} \label{lem:vp}
Fix $p>2$. Then, \eqref{eq:Psi-bound-again} holds for every $u\in[2,v_p]$.
\end{lemma}

\begin{proof}
If $p>6$, the conclusion is obvious as $v_p=2$ so we assume that $p\in(2,6]$. By Lemma \ref{lem:HY},
\begin{equation} \label{eq:use-HY}
\Psi_p(u) \leq \frac{\pi}{2\Gamma(1+\frac{1}{p})} \bigg( \frac{u}{u-1}\bigg)^{(u-1)(\frac12-\frac1p)} \leq \frac{\pi}{2\Gamma(1+\frac{1}{p})} \bigg( \frac{v_p}{v_p-1}\bigg)^{(v_p-1)(\frac12-\frac1p)},
\end{equation}
since the right-hand side is increasing in $u$. For $2< p<3$, the desired estimate becomes
\begin{equation}
-\big(\tfrac12+\tfrac1p\big) \log \big(\tfrac12+\tfrac1p\big) \leq \log\Big(\tfrac{2}{\pi}\Gamma\big(1+\tfrac1p\big)\Psi_p(\infty)\Big) \stackrel{\eqref{eq:derivation-small-coefficients}}{=} \tfrac12 \Big(\log \tfrac6\pi + 3 \log \Gamma\big(1+\tfrac1p\big) - \log\Gamma\big(1+\tfrac3p\big)\Big).
\end{equation}
Denoting by $a=\tfrac1p\in(\frac13,\frac12]$, the difference of the two sides
\begin{equation}
\delta(a) \eqdef -\big(\tfrac12+a\big) \log \big(\tfrac12+a\big)  - \tfrac12 \Big(\log \tfrac6\pi + 3 \log \Gamma\big(1+a\big) - \log\Gamma\big(1+3a\big)\Big)
\end{equation}
has $\delta(\frac12)=0$ and satisfies, by the standard formula \cite[6.3.16]{AS64},
\begin{equation}
\delta'(a) = -\log\big(\tfrac12+a\big) -1 + \frac12\sum_{k=1}^\infty \Big( \frac{3}{k+a} - \frac{3}{k+3a} \Big).
\end{equation}
In particular, for $a=\frac12$ the sum is telescoping and $\delta'(\frac12)=0$. Differentiating again,
\begin{equation}
\delta''(a) = -\frac{1}{\frac12+a}-\frac12 \Big( \sum_{k=1}^\infty \frac{3}{(k+a)^2} - \sum_{k=1}^\infty \frac{9}{(k+3a)^2} \Big)
\end{equation}
and
\begin{equation*}
\delta'''(a) = \frac{1}{(\frac12+a)^2}+  \sum_{k=1}^\infty \frac{3}{(k+a)^3} - \sum_{k=1}^\infty \frac{27}{(k+3a)^3} \leq \frac{36}{25} + \sum_{k=1}^\infty \frac{3}{(k+\frac13)^3} - \sum_{k=1}^\infty \frac{27}{(k+\frac32)^3} < -0.06.
\end{equation*}
As $\delta''(\frac13) \approx 0.05$ and $\delta''(\frac12) \approx -0.19$ and $\delta''$ is decreasing, there exists a number $a_0\in(\frac13,\frac12)$ such that $\delta''>0$ on $[\frac13,a_0)$ and $\delta''<0$ on $(a_0,\frac12]$. To conclude that $\delta\leq0$ observe that for $a\in[a_0,\frac12)$, Taylor's formula with a remainder gives
\begin{equation} \label{eq:a_0}
\delta(a) = \delta\big(\tfrac12\big)+\delta'\big(\tfrac12\big) \big(a-\tfrac12\big) + \int_a^{1/2} \delta''(u) (u-a)\,\diff u=  \int_a^{1/2} \delta''(u) (u-a)\,\diff u <0.
\end{equation}
On the other hand, $\delta$ is convex on the interval $[\tfrac13,a_0]$ with $\delta(\frac13)\approx-0.001$ and $\delta(a_0)<0$ by \eqref{eq:a_0}, so $\delta$ must be negative in the entire interval. This proves \eqref{eq:Psi-bound-again} in this range.

The range $3\leq p\leq 6$ is simpler.  Using again \eqref{eq:use-HY}, the desired estimate becomes
\begin{equation}
2\big(\tfrac12-\tfrac1p\big)\log\tfrac32  \leq \log\Big(\tfrac{2}{\pi}\Gamma\big(1+\tfrac1p\big)\Psi_p(\infty)\Big) \stackrel{\eqref{eq:derivation-small-coefficients}}{=} \tfrac12 \Big(\log \tfrac6\pi + 3 \log \Gamma\big(1+\tfrac1p\big) - \log\Gamma\big(1+\tfrac3p\big)\Big).
\end{equation}
By similar reasoning, the difference $\delta(a)$ of the two sides is convex in $a=\tfrac1p\in[\tfrac16,\tfrac13]$ and thus we conclude that $\delta(a)\leq \max\{\delta(1/6),\delta(1/3)\} <-0.001$ and  \eqref{eq:Psi-bound-again} follows.
\end{proof}

Having established \eqref{eq:Psi-bound-again} for $u\leq v_p$, we are left with proving the inequality for $u>v_p$. Consider the rescaled function
\begin{equation}
\widetilde\gamma_p(t) \eqdef \gamma_p\!\left(\frac{t}{\sqrt{c_p}}\right),
\end{equation}
where $c_p$ is given by \eqref{eq:cp}, and observe that \eqref{eq:Psi-bound-again} is equivalent to
\begin{equation} \label{eq:equi}
\begin{split}
\int_0^\infty\big|\widetilde\gamma_p(t)\big|^u\, \diff t & = \sqrt{c_p} \int_0^\infty \big|\gamma_p(s)\big|^u\, \diff s \stackrel{\eqref{eq:Psi-bound-again}}{\leq} \frac{\sqrt{c_p}}{\sqrt{u}} \max\big\{\Psi_p(2),\Psi_p(\infty)\big\} 
\\ & = \frac{\sqrt{\pi}}{2\sqrt{u}} \underbrace{\max\bigg\{ 2^{\frac12-\frac1p} \bigg( \frac{6}{\pi} \frac{\Gamma(1+\frac1p)^3}{\Gamma(1+\frac3p)} \bigg)^{-1/2},1\bigg\}}_{M_p}.
\end{split}
\end{equation}
Let $\lambda_p$ be the parameter given by \eqref{eq:lambdap}, and fix $u>v_p$. By Corollary \ref{cor:Qp}, if $|\widetilde{\gamma}_p(t)|>\lambda_p$, then $0<\widetilde\gamma_p(t) \leq (1-\theta_p t^2)^{1/\theta_p}$, where $\theta_p$ is given by \eqref{eq:thetap}. Hence, for every $t>0$, we have
\begin{equation} \label{eq:used-Qp-bb}
\begin{split}
\big|\widetilde\gamma_p(t)\big|^u  &= \lambda_p^{u-v_p} \big|\widetilde\gamma_p(t)\big|^{v_p} + \big|\widetilde\gamma_p(t)\big|^{v_p} \big( \big|\widetilde\gamma_p(t)\big|^{u-v_p} - \lambda_p^{u-v_p} \big)
\\ & \leq  \lambda_p^{u-v_p} \big|\widetilde\gamma_p(t)\big|^{v_p} + \big|\widetilde\gamma_p(t)\big|^{v_p} \big( \big|\widetilde\gamma_p(t)\big|^{u-v_p} - \lambda_p^{u-v_p} \big)_+
\\ & \leq  \lambda_p^{u-v_p} \big|\widetilde\gamma_p(t)\big|^{v_p} + \Big( \big(1-\theta_p t^2\big)_+^{u/\theta_p} - \lambda_p^{u-v_p} \big(1-\theta_p t^2\big)_+^{v/\theta_p} \Big)_+.
\end{split}
\end{equation}
Integrating over $(0,\infty)$, we get
\begin{equation*}
\begin{split}
\int_0^\infty \big|\widetilde\gamma_p(t)\big|^u \,\diff t \leq \lambda_p^{u-v_p} \int_0^\infty \big|\widetilde\gamma_p(t)\big|^{v_p}\,\diff t + \int_0^{\rho_p}  \Big( \big(1-\theta_p t^2\big)^{u/\theta_p} - \lambda_p^{u-v_p} \big(1-\theta_p t^2\big)^{v_p/\theta_p} \Big)\,\diff t,
\end{split}
\end{equation*}
where $\rho_p$ is determined by $(1-\theta_p\rho_p^2)^{1/\theta_p} = \lambda_p$. Using Lemma \ref{lem:vp} for $u=v_p$ and \eqref{eq:equi} yields
\begin{equation} \label{eq:used-vp}
\int_0^\infty \big|\widetilde\gamma_p(t)\big|^u \,\diff t \leq \frac{\sqrt{\pi}}{2\sqrt{v_p}} M_p \lambda_p^{u-v_p} + \int_0^{\rho_p}  \Big( \big(1-\theta_p t^2\big)^{u/\theta_p} - \lambda_p^{u-v_p} \big(1-\theta_p t^2\big)^{v_p/\theta_p} \Big)\,\diff t.
\end{equation}
The proof will be completed by the following proposition.

\begin{proposition} \label{prop:G-increasing}
Fix $p>2$. The function $G_p:[v_p,\infty)\to\R$ given by
\begin{equation} \label{eq:def-G}
\forall \ v\geq v_p,\qquad G_p(v) \eqdef \lambda_p^{-v}\bigg( \frac{\sqrt{\pi}}{2\sqrt{v}} - \int_0^{\rho_p} \big(1-\theta_pt^2\big)^{v/\theta_p}\,\diff t\bigg)
\end{equation}
is nondecreasing.
\end{proposition}

\begin{proof} [Proof of Theorem \ref{thm:Psi} using Proposition \ref{prop:G-increasing}]
Since $G_p$ is increasing, the same holds for the function
\begin{equation}
\Omega_p(v) \eqdef \lambda_p^{-v}\bigg( \frac{\sqrt{\pi}}{2\sqrt{v}}M_p - \int_0^{\rho_p} \big(1-\theta_pt^2\big)^{v/\theta_p}\,\diff t\bigg).
\end{equation}
Indeed, for $p\leq p_0$ we have $M_p=1$ and thus $\Omega_p \equiv G_p$, whereas for $p>p_0$ the derivative of $\Omega_p-G_p$ is a positive multiple of
\begin{equation}
\frac{\diff}{\diff v} \frac{\lambda_p^{-v}}{\sqrt v} = \frac{\lambda_p^{-v}}{\sqrt v}\Big(\log\tfrac1\lambda_p - \tfrac{1}{2v}\Big) \geq \frac{\lambda_p^{-v}}{\sqrt v}\Big(\log\tfrac83-\tfrac14\Big) > 0,
\end{equation}
where the first inequality uses that $\lambda_p= \frac{1}{3\Gamma(1+\frac1p)} \leq \frac{1}{3\Gamma(\frac43)} <\tfrac38$ for any $p\geq3$. Therefore, by \eqref{eq:used-vp},
\begin{equation*}
\begin{split}
\int_0^\infty \big|\widetilde\gamma_p(t)\big|^u \,\diff t  \leq  \lambda_p^u \Omega_p(v_p) +\! \int_0^{\rho_p} \big(1-\theta_p t^2\big)^{u/\theta_p}\,\diff t \leq \lambda_p^u \Omega_p(u) + \! \int_0^{\rho_p} \big(1-\theta_p t^2\big)^{u/\theta_p}\,\diff t
 = \frac{\sqrt{\pi}}{2\sqrt{u}} M_p
\end{split}
\end{equation*}
and this completes the proof.
\end{proof}

For parameters $\theta, v>0$ and $\lambda\in(0,1)$, consider the more general function 
\begin{equation} \label{eq:general-G}
G(\theta,\lambda,v) \eqdef \lambda^{-v} \bigg( \frac{\sqrt{\pi}}{2\sqrt{v}} - \int_0^{\rho(\theta,\lambda)} \big(1-\theta t^2\big)^{v/\theta}\,\diff t\bigg),
\end{equation}
where $\rho(\theta,\lambda)$ satisfies $(1-\theta\rho(\theta,\lambda)^2)^{1/\theta} = \lambda,$ or equivalently
\begin{equation} \label{eq:rho-theta}
    \rho(\theta,\lambda) = \sqrt{\frac{1-\lambda^\theta}{\theta}}.
\end{equation}
We start with the following simpler lemma.

\begin{lemma} \label{lem:G-convex}
Fix $\theta>0$ and $\lambda\in(0,1)$. The function $G(\theta,\lambda,\cdot)$ is convex on $(0,\infty)$.
\end{lemma}

\begin{proof}
By the Gaussian identity $\int_0^\infty e^{-vt^2}\,\diff t= \frac{\sqrt\pi}{2\sqrt v}$, we can write
\begin{equation} \label{eq:G-in-2}
G(\theta,\lambda,v) = \int_0^{\rho(\theta,\lambda)} \bigg[ \Big( \frac{e^{-t^2}}{\lambda} \Big)^v - \Big( \frac{(1-\theta t^2)^{1/\theta}}{\lambda} \Big)^v \bigg] \, \diff t + \int_{\rho(\theta,\lambda)}^\infty \Big( \frac{e^{-t^2}}{\lambda} \Big)^v  \, \diff t.
\end{equation}
Differentiating twice gives
\begin{equation}
\begin{split}
\frac{\partial^2 G}{\partial v^2}(\theta,\lambda,v)  = \int_0^{\rho(\theta,\lambda)} \bigg[ \Big( \frac{e^{-t^2}}{\lambda} \Big)^v \log^2\Big(\frac{e^{-t^2}}{\lambda}\Big) - \Big( &\frac{(1-\theta t^2)^{1/\theta}}{\lambda}\Big)^v  \log^2\Big(\frac{(1-\theta t^2)^{1/\theta}}{\lambda}\Big) \bigg] \, \diff t 
\\ & + \int_{\rho(\theta,\lambda)}^\infty \Big( \frac{e^{-t^2}}{\lambda}\Big)^v \log^2\Big(\frac{e^{-t^2}}{\lambda}\Big)   \, \diff t.
\end{split}
\end{equation}
By the definition of $\rho$, for $0\leq t\leq\rho(\theta,\lambda)$ we have
\begin{equation} \label{eq:Berno}
\frac{e^{-t^2}}{\lambda} \geq \frac{(1-\theta t^2)^{1/\theta}}{\lambda} \geq 1
\end{equation}
and thus $\frac{\partial^2G}{\partial v^2} \geq0$ since the function $x\mapsto x^v \log^2x$ is increasing on $[1,\infty)$.
\end{proof}

Given the convexity of $G$, Proposition \ref{prop:G-increasing} will follow once we prove that $G_p'(v_p)\geq0$. For parameters $\theta, v>0$ and $\lambda\in(0,1)$, consider the function 
\begin{equation} \label{eq:H-function}
    \begin{split} 
H(\theta,\lambda,v) \eqdef \int_0^{\rho(\theta,\lambda)} \Big(\frac{(1-\theta t^2)^{1/\theta}}{\lambda}\Big)^v \log\Big(\frac{(1-\theta t^2)^{1/\theta}}{\lambda} \Big)\, \diff t.
\end{split}
\end{equation}
Differentiating under the integral in the definition \eqref{eq:def-G} of $G_p$ gives
\begin{equation} \label{eq:G'-formula}
    \frac{\partial G}{\partial v}(\theta,\lambda,v) = \frac{\sqrt{\pi}}{2\sqrt{v}\lambda^v}\big(\log(\tfrac{1}{\lambda}) - \tfrac{1}{2v}\big) - H(\theta,\lambda, v).
\end{equation}

We will use the following lemma.

\begin{lemma} \label{lem:H-lemma}
 For fixed $\lambda\in(0,1)$ and $v>0$, the map $\theta\mapsto H(\theta,\lambda,v)$ is decreasing on $(0,\infty)$. 
\end{lemma}

\begin{proof}
First, observe that $\theta\mapsto\rho(\theta,\lambda)$ is decreasing for $\theta>0$, since the function $\theta\mapsto \lambda^\theta$ is convex and 
\begin{equation}
    \rho(\theta,\lambda)^2 = \frac{1-\lambda^\theta}{\theta} = -\frac{\lambda^\theta - \lambda^0}{\theta}.
\end{equation}
Moreover, fix $t>0$ and consider the function $\alpha_t(\theta) = (1-\theta t^2)^{1/\theta}$ defined for $\theta\in[0,1/t^2]$. Then, we can expand
\begin{equation}
    \log\alpha_t(\theta) = \frac{\log(1-\theta t^2)}{\theta} = -\sum_{n=1}^\infty \frac{t^{2n} \theta^{n-1}}{n},
\end{equation}
which readily implies that $\alpha_t$ is decreasing as well. Finally, consider the function $\psi_{v}(w) = w^v \log w$ defined for $w>0$, which is nonnegative on $[1,\infty)$ and is increasing on $[ e^{-1/v},\infty)$. Consider now $0<\theta_1\leq\theta_2$ and write the function as
\begin{equation}
    H(\theta_2,\lambda,v) = \int_0^{\rho(\theta_2,\lambda)} \psi_{v}\Big(\frac{\alpha_t(\theta_2)}{\lambda}\Big) \,\diff t.
\end{equation}
Fix $t\in[0,\rho(\theta_2,\lambda)] \subseteq [0,\rho(\theta_1,\lambda)]$ since $\rho(\cdot,\lambda)$ is decreasing. Since
\begin{equation}
    t\leq \rho(\theta_2,\lambda) = \sqrt{\frac{1-\lambda^{\theta_2}}{\theta_2}} \leq \frac{1}{\sqrt{\theta_2}},
\end{equation}
we have $\theta_1\leq \theta_2\leq t^{-2}$ and the fact that $\alpha_t$ is decreasing yields $\alpha_t(\theta_2) \leq \alpha_t(\theta_1)$. Moreover, by the definition \eqref{eq:rho-theta} of $\rho(\theta,\lambda)$, the condition $t\in[0,\rho(\theta_2,\lambda)]$ shows that
\begin{equation}
    \alpha_t(\theta_1)\geq \alpha_t(\theta_2) = (1-\theta_2 t^2)^{1/\theta_2} \geq \lambda \geq \lambda e^{-1/v}.
\end{equation}
Thus, the monotonicity of $\psi_v$ on $[e^{-1/v},\infty)$ yields $\psi_v(\alpha_t(\theta_2)/\lambda) \leq \psi_v(\alpha_t(\theta_1)/\lambda)$. Integrating, we get
\begin{equation}
\begin{split}
H(\theta_2,\lambda,v) = \int_0^{\rho(\theta_2,\lambda)} \psi_{v}\Big(\frac{\alpha_t(\theta_2)}{\lambda}\Big) \,\diff t 
 & \leq  \int_0^{\rho(\theta_2,\lambda)} \psi_{v}\Big(\frac{\alpha_t(\theta_1)}{\lambda}\Big) \,\diff t 
\\ & \leq  \int_0^{\rho(\theta_1,\lambda)} \psi_{v}\Big(\frac{\alpha_t(\theta_1)}{\lambda}\Big) \,\diff t = H(\theta_1,\lambda, v),
\end{split}
\end{equation}
where the second inequality follows since $\rho(\theta_2,\lambda)\leq \rho(\theta_1,\lambda)$ and because for $t\in[\rho(\theta_2,\lambda),\rho(\theta_1,\lambda)]$,
\begin{equation*}
    \frac{\alpha_t(\theta_1)}{\lambda} \geq 1 \quad \Longrightarrow \quad \psi_v\Big(\frac{\alpha_t(\theta_1)}{\lambda}\Big) \geq0. \hfill\qedhere
\end{equation*}
\end{proof}

We are now ready to complete the proof of the Proposition.

\begin{proof} [Proof of Proposition \ref{prop:G-increasing}]
We shall distinguish two cases.

\smallskip

\noindent$\bullet$ \emph{Suppose that $p\geq3$.} In this range we have $\lambda_p=\frac{1}{3\Gamma(1+\frac1p)}\leq\frac{1}{3\Gamma(\frac43)}<\frac38$ and we consider
\begin{equation} \label{eq:Gtilde}
\forall\ v\geq2,\qquad \widetilde{G}_p(v)\stackrel{\mathrm{def}}{=}\big(\tfrac{8\lambda_p}{3}\big)^vG_p(v).
\end{equation}
The representation \eqref{eq:G-in-2} along with inequality \eqref{eq:Berno} imply that $G_p$ is positive and thus so is $\widetilde{G}_p$. Therefore, as $\log\frac{3}{8\lambda_p}>0$, differentiating the identity $G_p(v)=(\frac{3}{8\lambda_p})^v\widetilde{G}_p(v)$ gives
\begin{equation} \label{eq:Gprime}
G_p'(v)=\big(\tfrac{3}{8\lambda_p}\big)^v\Big(\widetilde{G}_p'(v)+\log\big(\tfrac{3}{8\lambda_p}\big)\widetilde{G}_p(v)\Big)\geq\big(\tfrac{3}{8\lambda_p}\big)^v\widetilde{G}_p'(v).
\end{equation}
Since $G_p$ is convex by Lemma \ref{lem:G-convex}, it suffices to prove that $\widetilde{G}_p'(v_p)\geq0$. By \eqref{eq:rho-theta}, the function $\rho(\theta_p,\cdot)$ is decreasing on $(0,1)$, so $\rho(\theta_p,\frac38)<\rho(\theta_p,\lambda_p)$ and the definitions \eqref{eq:def-G} and \eqref{eq:general-G} give
\begin{equation} \label{eq:Gtilde-split}
\widetilde{G}_p(v)=G\big(\theta_p,\tfrac38,v\big)-\big(\tfrac83\big)^v\int_{\rho(\theta_p,\frac38)}^{\rho(\theta_p,\lambda_p)}\big(1-\theta_pt^2\big)^{v/\theta_p}\,\diff t.
\end{equation}
Differentiating \eqref{eq:Gtilde-split} and using (108) with $\lambda=\frac38$, we get
\begin{equation} \label{eq:Gtilde-der}
\begin{split}\widetilde{G}_p'(v)=\big(\tfrac83\big)^v\frac{\sqrt{\pi}}{2\sqrt{v}}\Big(\log\tfrac83&-\frac{1}{2v}\Big)-H\big(\theta_p,\tfrac38,v\big)\\&-\int_{\rho(\theta_p,\frac38)}^{\rho(\theta_p,\lambda_p)}\Big(\tfrac83\big(1-\theta_pt^2\big)^{1/\theta_p}\Big)^v\log\Big(\tfrac83\big(1-\theta_pt^2\big)^{1/\theta_p}\Big)\,\diff t.\end{split}
\end{equation}
By the definition \eqref{eq:rho-theta} of $\rho(\theta,\lambda)$, for every $t\in[\rho(\theta_p,\frac38),\rho(\theta_p,\lambda_p)]$ we have $\lambda_p\leq(1-\theta_pt^2)^{1/\theta_p}\leq\frac38$, so the logarithm in the last integrand is nonpositive and \eqref{eq:Gtilde-der} yields
\begin{equation} \label{eq:Gtilde-lower}
\widetilde{G}_p'(v)\geq\big(\tfrac83\big)^v\frac{\sqrt{\pi}}{2\sqrt{v}}\Big(\log\tfrac83-\frac{1}{2v}\Big)-H\big(\theta_p,\tfrac38,v\big).
\end{equation}
If $p\geq6$, then $\theta_p\geq\theta_6=\frac13$ since $p\mapsto\theta_p$ is increasing by Lemma \ref{lem:thetap}, thus \eqref{eq:Gtilde-lower} combined with Lemma \ref{lem:H-lemma} gives
\begin{equation}
\widetilde{G}_p'(2)\geq\big(\tfrac83\big)^2\frac{\sqrt{\pi}}{2\sqrt{2}}\Big(\log\tfrac83-\frac14\Big)-H\big(\tfrac13,\tfrac38,2\big)\approx0.011>0.
\end{equation}
If $3\leq p\leq6$, then we have $\theta_p\geq\theta_3=1-\frac{4\pi}{9\sqrt{3}}$, thus \eqref{eq:Gtilde-lower} combined with Lemma \ref{lem:H-lemma} gives
\begin{equation}
\widetilde{G}_p'(3)\geq\big(\tfrac83\big)^3\frac{\sqrt{\pi}}{2\sqrt{3}}\Big(\log\tfrac83-\frac16\Big)-H\big(1-\tfrac{4\pi}{9\sqrt{3}},\tfrac38,3\big)\approx0.076>0.
\end{equation}

\noindent $\bullet$ {\it Suppose that $p\in(2,3)$.} This case is more intricate. First observe that since $a\mapsto\theta_{1/a}$ is concave on $(0,\frac12)$, for $p\in[2,6]$ we have
\begin{equation}
\theta_p \geq \big(\tfrac3p-\tfrac12\big) \theta_2 + \big(\tfrac32-\tfrac3p\big) \theta_6 = \frac{p-2}{2p},
\end{equation}
where we used again that $\theta_6 = \frac13$. Therefore, by Lemma \ref{lem:H-lemma}
\begin{equation}
G_p'(v) \stackrel{\eqref{eq:G'-formula}}{=} \frac{\sqrt{\pi}2^v}{2\sqrt{v}}\big(\log2 - \tfrac{1}{2v}\big) - H(\theta_p,\lambda_p, v) \geq \frac{\sqrt{\pi}2^v}{2\sqrt{v}}\big(\log2 - \tfrac{1}{2v}\big) - H\big(\tfrac{p-2}{2p},\tfrac12, v\big).
\end{equation}
In particular, plugging $v=v_p=\frac{2p}{p-2}$, this becomes
\begin{equation} \label{eq:G'-Theta}
\begin{split}
\frac{G_p'(v_p)}{2^{v_p}} & \geq \frac{\sqrt{\pi}}{2\sqrt{v_p}}\big(\log2 - \tfrac{1}{2v_p}\big) -  2^{-v_p} H\big(\tfrac{1}{v_p},\tfrac12, v_p\big)
\\ & =  \frac{\sqrt{\pi}}{2\sqrt{v_p}}\big(\log2 - \tfrac{1}{2v_p}\big) - \int_0^{\sqrt{v_p(1-2^{-1/v_p})}} \Big(1-\frac{t^2}{v_p}\Big)^{v_p^2} \log\Big(2 \Big(1-\frac{t^2}{v_p}\Big)^{v_p} \Big) \, \diff t,
\end{split}
\end{equation}
Take now an arbitrary $v>6$. Then,
\begin{equation} \label{eq:lower-theta}
\begin{split}
\Theta(v) & \eqdef \frac{\sqrt{\pi}}{2\sqrt{v}} \big(\log2 - \tfrac{1}{2v}\big)  -  2^{-v} H\big(\tfrac{1}{v},\tfrac12, v\big)
\\ & =  \int_0^\infty e^{-vt^2}\big(\log2 - t^2\big)\,\diff t  - \int_0^{\sqrt{v(1-2^{-1/v})}} \Big(1-\frac{t^2}{v}\Big)^{v^2} \log\Big(2 \Big(1-\frac{t^2}{v}\Big)^{v} \Big) \, \diff t
\\ & \geq \int_0^\infty e^{-vt^2}\big(\log2 - t^2\big)\,\diff t- \int_0^{\sqrt{v(1-2^{-1/v})}} e^{-vt^2-t^4/2} \big(\log2 - t^2 \big) \, \diff t
\end{split}
\end{equation}
by two applications of the inequality $(1-\frac{t^2}{v})^v \leq e^{-t^2-t^4/2v}\leq e^{-t^2}$. Moreover, observe that $v(1-2^{-1/v}) = v(1-e^{-\log2/v}) \leq \log 2$. Consider the 3 consecutive intervals $I_1 =[0,\sqrt{v(1-2^{-1/v})}]$, $I_2 = [v(1-2^{-1/v}),\sqrt{\log2}]$ and $I_3 = [\sqrt{\log2},\infty)$. The lower bound \eqref{eq:lower-theta} for $\Theta$ is a sum of 3 terms, the first of which satisfies
\begin{equation}
\Sigma_1 = \int_{I_1} e^{-vt^2}\big(1-e^{-t^4/2}\big)\big(\log2 - t^2\big)\, \diff t \geq \tfrac{7}{16} \int_{I_1} t^4e^{-vt^2}\big(\log2 - t^2\big)\, \diff t,
\end{equation}
since $1-e^{-y} \geq y-y^2/2>\frac{7y}{8}$ for $y\in[0,\frac14]$. The second satisfies
\begin{equation}
\Sigma_2 =\int_{I_2} e^{-vt^2}\big(\log2 - t^2\big)\,\diff t \geq \tfrac{7}{16} \int_{I_2} t^4e^{-vt^2}\big(\log2 - t^2\big)\, \diff t 
\end{equation}
since $\frac{7t^4}{16}  <1$ for $t\leq \sqrt{\log2}$. Finally, the third integral can be written as
\begin{equation}
\begin{split}
\Sigma_3  & = \tfrac{7}{16}\int_{I_3}t^4 e^{-vt^2}\big(\log2 - t^2\big)\,\diff t - \int_{\sqrt{\log2}}^\infty e^{-vt^2} \big( t^2-\log 2\big) \Big(1-\frac{7t^4}{16}\Big) \, \diff t
\\ & \geq \tfrac{7}{16}\int_{I_3}t^4 e^{-vt^2}\big(\log2 - t^2\big)\,\diff t - \int_{\sqrt{\log2}}^\infty e^{-vt^2} \big( t^2-\log 2\big) \Big(1+\frac{7t^4}{16}\Big) \, \diff t
\end{split}
\end{equation}
Using the inequality $1+\frac{7t^4}{16} \leq 1+\frac{t^4}{2} \leq e^{t^2}$ and substituting $t=s+\sqrt{\log 2}$, we can estimate the last integral as
\begin{equation}
\begin{split}
\int_{\sqrt{\log2}}^\infty e^{-vt^2}& \big( t^2-\log 2\big) \Big(1+\frac{7t^4}{16}\Big) \, \diff t\leq \int_{\sqrt{\log2}}^\infty e^{-(v-1)t^2} \big( t^2-\log 2\big) \, \diff t
\\ & = 2^{-v+1} \int_0^\infty e^{-(v-1)(s^2+2\sqrt{\log2}s)} \big(s^2 + 2\sqrt{\log 2}s\big)\,\diff s
\\ & \leq 2^{-v+1} \int_0^\infty e^{-2\sqrt{\log 2}(v-1)s} \big(s^2 + 2\sqrt{\log 2}s\big)\,\diff s
\\ & = 2^{-v}\Big(\frac{1}{(v-1)^2\sqrt{\log2}}+\frac{1}{2(v-1)^3(\log2)^{3/2}}\Big).
\end{split}
\end{equation}
Combining all the above, we finally conclude that
\begin{equation} \label{eq:Theta-final}
\begin{split}
2^v\Theta(v) & \geq \tfrac{7}{16} 2^v \int_{0}^\infty t^4 e^{-vt^2}\big(\log2 - t^2\big)\,\diff t - \frac{1}{(v-1)^2\sqrt{\log2}}-\frac{1}{2(v-1)^3(\log2)^{3/2}}
\\ & = \frac{21\sqrt{\pi}}{128}\Big(\log2-\frac{5}{2v}\Big)\frac{2^v}{v^{5/2}} - \frac{1}{(v-1)^2\sqrt{\log2}} - \frac{1}{2(v-1)^3(\log2)^{3/2}}.
\end{split}
\end{equation}
The first term is increasing for $v\geq \frac{5}{2\log2} \approx 3.61$ while the other two are decreasing. Bounding below by the value of all these terms at $v=6$ and combining \eqref{eq:G'-Theta}, \eqref{eq:Theta-final}, we finally get
\begin{equation}
G_p'(v_p) \geq 2^{v_p}\Theta(v_p) \geq \frac{21\sqrt{\pi}}{128}\Big(\log2-\frac{5}{12}\Big)\frac{2^6}{6^{5/2}} - \frac{1}{25\sqrt{\log2}} - \frac{1}{250(\log2)^{3/2}} > 0.003
\end{equation}
and this completes the proof of the Proposition.
\end{proof}

As mentioned earlier, we shall further need a deficit term of inequality \eqref{eq:Psi-bound-again} for $p$ far from 2.

\begin{theorem} \label{thm:Psi-deficit}
Fix $p>26$. Then, for every $u\geq2$, we have
\begin{equation} \label{eq:Psi-bound-deficit}
1-\frac{\Psi_p(u)}{\max\{ \Psi_p(2), \Psi_p(\infty)\}} \geq \frac{3}{100} \frac{u-2}{u^2}.
\end{equation}
\end{theorem}

\begin{proof}
The proof requires a more careful analysis of the argument used for Theorem \ref{thm:Psi} in the large $p$ regime. Since $p>26$, we have $\lambda_p=\frac{1}{3\Gamma(1+\frac1p)} < \frac{1}{3\Gamma(\frac{27}{26})}<\frac{7}{20} \eqdef \lambda_*$. Incorporating this bound in \eqref{eq:used-Qp-bb}, integrating and recalling the definition \eqref{eq:Psip}, we obtain
\begin{equation} \label{eq:TBV1}
\frac{\sqrt{c_p}\Psi_p(u)}{\sqrt u} \leq \int_0^{\rho(\theta_p,\lambda_*)}\big(1-\theta_pt^2\big)^{u/\theta_p}\, \diff t +\lambda_*^{u-2} \bigg( \frac{\sqrt{c_p}\Psi_p(2)}{\sqrt2} -\int_0^{\rho(\theta_p,\lambda_*)}\big(1-\theta_pt^2\big)^{2/\theta_p}\, \diff t \bigg).
\end{equation}
Equivalently, since $\sqrt{c_p}\Psi_p(\infty) = \frac{\sqrt{\pi}}{2}$, this can be written in terms of the function $G$ of \eqref{eq:general-G} as
\begin{equation} \label{eq:main-for-def}
\begin{alignedat}{2}
\sqrt{c_p}&\big( \max\big\{\Psi_p(2),\Psi_p(\infty)\big\} &&- \Psi_p(u)\big)
\\ & \mathrlap{\stackrel{\eqref{eq:TBV1}}{\geq} \bigg(\frac{\sqrt{\pi}}{2} - \sqrt{u} \int_0^{\rho(\theta_p,\lambda_*)}\big(1-\theta_pt^2\big)^{u/\theta_p} \, \diff t\bigg) - \sqrt{u} \lambda_*^{u-2}\bigg( \frac{\sqrt{\pi}}{2\sqrt{2}} -\int_0^{\rho(\theta_p,\lambda_*)}\big(1-\theta_pt^2\big)^{2/\theta_p}\, \diff t \bigg)}
\\ & && \hspace{-6mm} +\sqrt{c_p} \big( \max\big\{\Psi_p(2), \Psi_p(\infty) \big\} - \Psi_p(\infty) \big) - \frac{\sqrt{c_pu}}{\sqrt2}\lambda_*^{u-2} \big(\Psi_p(2) - \Psi_p(\infty)\big)
\\ & \mathrlap{\stackrel{\eqref{eq:general-G}}{=} \sqrt{u}\lambda_*^u\big(G(\theta_p,\lambda_*,u)-G(\theta_p,\lambda_*,2)\big)}
\\ & && \hspace{-6mm} +\sqrt{c_p} \big( \max\big\{\Psi_p(2), \Psi_p(\infty) \big\} - \Psi_p(\infty) \big) -\frac{\sqrt{c_pu}}{\sqrt2}\lambda_*^{u-2} \big(\Psi_p(2) - \Psi_p(\infty)\big).
\end{alignedat}
\end{equation}

\noindent $\bullet$ \emph{Suppose that $2\leq u\leq 4$.} If $p\geq p_0$, then $\Psi_p(2)\geq\Psi_p(\infty)$ and the sum of the last two terms is a positive multiple of
\begin{equation}
1 - \frac{\sqrt u}{\sqrt2} \lambda_*^{u-2} = 1 - \frac{\sqrt u}{\sqrt2} \Big(\frac{7}{20}\Big)^{u-2} \geq0
\end{equation}
for every $u\geq2$. In the opposite case, $\Psi_p(2)<\Psi_p(\infty)$ so the first term vanishes and what remains is evidently nonnegative. Hence, in either case, we conclude that
\begin{equation}
\begin{split}
\sqrt{c_p}\big( \max\big\{\Psi_p(2),\Psi_p(\infty)\big\} - \Psi_p(u)\big) & \geq \sqrt{u}\Big(\frac{7}{20}\Big)^u \big(G\big(\theta_p,\tfrac{7}{20},u\big)-G\big(\theta_p,\tfrac{7}{20},2\big)\big) 
\\ & \geq\sqrt{u}(u-2)\Big(\frac{7}{20}\Big)^u \  \frac{\partial G}{\partial v}\big(\theta_p,\tfrac{7}{20},2\big),
\end{split}
\end{equation}
where we used that $G(\theta_p,\frac{7}{20},\cdot)$ is convex on $[2,\infty)$ by Lemma \ref{lem:G-convex}. By \eqref{eq:G'-formula} and Lemma \ref{lem:H-lemma},
\begin{equation}
\begin{split}
\frac{\partial G}{\partial v}\big(\theta_p,\tfrac{7}{20},2\big)  & =   \frac{200\sqrt{\pi}}{49\sqrt{2}}\big(\log(\tfrac{20}{7}) - \tfrac{1}{4}\big) - H\big(\theta_p,\tfrac{7}{20}, 2\big)
\\ & \geq \frac{200\sqrt{\pi}}{49\sqrt{2}}\big(\log(\tfrac{20}{7}) - \tfrac{1}{4}\big) - H\big(\tfrac{39}{100},\tfrac{7}{20}, 2\big) > \frac{1}{16}
\end{split}
\end{equation}
since $\theta_p \geq \theta_{26} > \frac{39}{100}$ and by Lemma \ref{lem:thetap}. Combining the above with the identity 
\begin{equation}
\sqrt{c_p}\max\big\{\Psi_p(2),\Psi_p(\infty)\big\} = \frac{\sqrt{\pi}M_p}{2}
\end{equation}
gives, for every $u\in[2,4]$,
\begin{equation}
1-\frac{\Psi_p(u)}{\max\{\Psi_p(2),\Psi_p(\infty)\}} \geq \frac{u^{5/2}}{8\sqrt{\pi}M_p}\Big(\frac{7}{20}\Big)^u \cdot \frac{u-2}{u^2} \geq \frac{\sqrt{3} u^{5/2}}{8\pi}\Big(\frac{7}{20}\Big)^u \cdot \frac{u-2}{u^2} \geq \frac{3}{100}\frac{u-2}{u^2},
\end{equation}
where we also used $M_p\leq\sqrt{\pi/3}$. To see this, observe that for $a\in[0,\frac12]$, the logarithm of the first term in the maximum defining $M_p$ in \eqref{eq:equi}, with $a=\frac1p$, equals
\begin{equation}
\big(\tfrac12-a\big)\log2-\tfrac12\Big(\log\tfrac6\pi+3\log\Gamma(1+a)-\log\Gamma(1+3a)\Big),
\end{equation}
which is a convex function of $a$ by the classical formula \cite[6.4.10]{AS64} for the derivatives of the digamma function. Hence it is at most the maximum of its values at $a=\frac12$ and $a=0$, which are $0$ and $\frac12\log\frac\pi3$ respectively. 

\smallskip

\noindent $\bullet$ \emph{Suppose that $u\geq4$.} Recalling the definition of $M_p\geq1$ from \eqref{eq:equi}, we can rewrite \eqref{eq:main-for-def} as
\begin{equation} \label{eq:Psi-def-2}
\begin{split}
1-&\frac{\Psi_p(u)}{\max\{\Psi_p(2),\Psi_p(\infty)\}}  \geq \frac{2\sqrt u}{\sqrt{\pi} M_p} \lambda_*^u G(\theta_p,\lambda_*,u)  +\frac{M_p-1}{M_p}
\\ &  +\frac{2\sqrt{u}}{\sqrt{\pi}M_p} \lambda_*^{u-2} \int_0^{\rho(\theta_p,\lambda_*)}\big(1-\theta_pt^2\big)^{2/\theta_p}\, \diff t  - \frac{\sqrt{u}}{\sqrt 2} \lambda_*^{u-2} \frac{\Psi_p(2)}{\max\{\Psi_p(2),\Psi_p(\infty)\}}.
\end{split}
\end{equation}
Since $\rho(\theta_p,\lambda_*)<\frac{1}{\sqrt{\theta_p}}$ by \eqref{eq:rho-theta},
\begin{equation}
\lambda_*^uG(\theta_p,\lambda_*,u) \stackrel{\eqref{eq:general-G}}{\geq} \frac{\sqrt\pi}{2\sqrt u} - \int_0^{\frac{1}{\sqrt{\theta_p}}} (1-\theta_pt^2)^{u/\theta_p}\, \diff t = \frac{\sqrt\pi}{2\sqrt u} \bigg( 1-\sqrt{\frac{u}{\theta_p}}\,
\frac{\Gamma(1+\frac{u}{\theta_p})}{\Gamma(\frac32+\frac{u}{\theta_p})}\bigg)
\end{equation}
by the usual relation between $B$ and $\Gamma$ functions. By an inequality of \cite{Ker83},
\begin{equation}
\sqrt{x}\, \frac{\Gamma(x+1)}{\Gamma(x+\frac32)} \leq \Big(1+\frac{3}{4x}\Big)^{-\frac12} \leq 1-\frac{1}{3x},
\end{equation}
where the last step holds for $x\geq10$. Since $\theta_p< \frac25$, we can apply this to $x=\frac{u}{\theta_p}\geq10$ to obtain
\begin{equation}
\lambda_*^uG(\theta_p,\lambda_*,u) \geq \frac{\sqrt\pi}{2\sqrt u}\frac{\theta_p}{3u}.
\end{equation}
Hence, the first two terms in the right-hand side of \eqref{eq:Psi-def-2} satisfy
\begin{equation}
\frac{2\sqrt u}{\sqrt{\pi} M_p} \lambda_*^u G(\theta_p,\lambda_*,u) +\frac{M_p-1}{M_p} \geq \frac{\theta_p}{3uM_p}  + \frac{M_p-1}{M_p}\geq\frac{\theta_p}{3u} \geq \frac{13}{100 u},
\end{equation}
where we used that $\frac{\theta_p}{3u} \leq1$.

To control the remaining two terms in \eqref{eq:Psi-def-2}, observe that both $\rho(\theta,\lambda_*)$ and $(1-\theta t^2)^{1/\theta}$ are decreasing functions of $\theta$ and $\theta_p<\frac25$, we therefore have
\begin{equation}
\int_0^{\rho(\theta_p,\lambda_*)} (1-\theta_pt^2)^{2/\theta_p}\,\diff t \geq \int_0^{\rho(\frac25,\frac{7}{20})} \big(1-\tfrac25t^2\big)^5\,\diff t > \frac{27}{50}.
\end{equation}
Now, if $\Psi_p(2)\leq\Psi_p(\infty)$, then
\begin{equation}
\frac{\Psi_p(2)}{\max\{\Psi_p(2),\Psi_p(\infty)\}}-\frac{2\sqrt2}{\sqrt\pi M_p}\int_0^{\rho(\theta_p,\lambda_*)} (1-\theta_pt^2)^{2/\theta_p}\,\diff t \leq 1-\frac{27\sqrt2}{25\sqrt\pi} < \frac{4}{25},
\end{equation}
whereas if $\Psi_p(2)>\Psi_p(\infty)$, then the first term is equal to $1$ and $M_p\leq\sqrt{\pi/3}$, so the left-hand side is at most $1-\frac{27\sqrt6}{25\pi}<\frac{4}{25}$. 

Plugging these estimates in \eqref{eq:Psi-def-2}, we get that for every $u\geq4$,
\begin{equation}
1-\frac{\Psi_p(u)}{\max\{\Psi_p(2),\Psi_p(\infty)\}} \geq \frac{13}{100u}-\frac{4}{25}\sqrt{\frac u2}\,\lambda_*^{u-2}.
\end{equation}
Since $\frac{13}{100u}-\frac{3(u-2)}{100u^2}=\frac{5u+3}{50u^2}$,  it now suffices to prove that
\begin{equation}
4\sqrt2\,\frac{u^{5/2}}{5u+3}\,\Big(\frac{7}{20}\Big)^{u-2} \leq 1.
\end{equation}
The logarithmic derivative of the left-hand side equals $\frac{15(u+1)}{2u(5u+3)}-\log\frac{20}{7}$, which is negative for $u\geq4$ as the first term is decreasing in $u$ and equal to $\frac{75}{184}<1<\log\frac{20}{7}$ at $u=4$. Hence, it suffices to check the inequality for $u=4$, where it reads $\frac{392\sqrt2}{575}<1$. This concludes the proof.
\end{proof}

\begin{remark}
The behavior of the Ball--Haagerup special function is much simpler for $p\in[1,2]$. Indeed, Koldobsky has proven \cite{Kol98} that for $p$ in this range, the function $s\mapsto \gamma_p(\sqrt{s})$ is completely monotonic on $[0,\infty)$ (see also \cite[p.~223]{BN02} for a different proof). Therefore, by Bernstein's theorem (see, e.g., \cite{Fel71}), there exists a measure $\mu$ on $(0,\infty)$ such that for every $s\in\R$,
\begin{equation}
\gamma_p(s) = \int_0^\infty e^{-t s^2} \, \diff\mu(t).
\end{equation}
Therefore, for $u\geq1$, we have
\begin{equation}
\Psi_p(u) = \int_0^\infty \big|\gamma_p(s/\sqrt{u})\big|^u \, \diff s  = \int_0^\infty \big\| t\mapsto e^{-t s^2}\big\|_{L_{1/u}(\mu)} \, \diff s.
\end{equation}
By H\"older's inequality, for any function $h$ on $(0,\infty)$, the function $u\mapsto \|h\|_{L_{1/u}(\mu)}$ is log-convex  and thus $\Psi_p$ is also log-convex as an integral of log-convex functions. In particular, the maximal value of $\Psi_p$ in any subinterval of $[1,\infty]$ is attained at the endpoints.
\end{remark}

\begin{remark}
In the concluding Remark of \cite{Kon25}, K\"onig points out that in order to prove the estimate \eqref{eq:Psi-bound-again} in the range $2<p<20$ via his distribution function argument, a more delicate approximation of $\gamma_p$ by damped trigonometric polynomials would be needed. In contrast, the bootstrapping argument of this section circumvents this need, relying only on the polynomial envelope of Proposition \ref{prop:Qp}, and even yields a stability version for sufficiently large $p$.
\end{remark}


\section{Proof of Theorem \ref{thm:main}} 

We are finally ready to complete the proof of the main result of this paper. In the sequel, we shall always tacitly assume that the unit vector $a$ satisfies $a_1\geq \cdots \geq a_n\geq0$. As explained in \eqref{eq:derivation-small-coefficients}, if $a_1 \leq \frac{1}{\sqrt{2}}$, the conclusion of Theorem \ref{thm:main} follows readily from Lemma \ref{lem:holder-Psi} and Theorem \ref{thm:Psi}. For the rest of the proof, we shall assume that $a_1> \frac{1}{\sqrt{2}}$. Denote by
\begin{equation}
T_p \eqdef 2^{\frac12-\frac1p} \qquad \mbox{and} \qquad G_p \eqdef \bigg( \frac{6}{\pi} \frac{\Gamma(1+\frac1p)^3}{\Gamma(1+\frac3p)} \bigg)^{1/2}
\end{equation}
so that the desired bound \eqref{eq:final-goal} becomes
\begin{equation} \label{eq:final-goal-2} 
A_{n,p}(a) \leq \max\{T_p, G_p\}.
\end{equation}
\medskip


\noindent {\bf Case 1. Suppose that ${\bf 2<p<4}$.} We shall use the following numerical lemma.

\begin{lemma}	\label{lem:large-coordinate-1}
For every $2<p < 4$, we have
\begin{equation} \label{eq:physical-scalar}
\exp\left\{\big(\tfrac14-\tfrac{1}{2p}\big)(1+\log2)\right\}\le G_p.
\end{equation}
\end{lemma}

The proof is deferred to Section \ref{sec:technical}.

\begin{proof} [Proof of Case 1.]
Consider the function $\phi:[0,1)\to\R_+$ given by
\begin{equation}
\forall \ t\in[0,1),\qquad	\phi(t):=-(1-t)\log(1-t).
\end{equation}
Barthe's estimate \eqref{eq:barthe} for the section function yields
\begin{equation}\label{eq:BL}
\log A_{n,p}(a)\le\big(\tfrac12-\tfrac{1}{p}\big) \sum_{i=1}^n \phi(a_i^2).
\end{equation}
Write
\begin{equation}
r \eqdef 1-a_1^2 = \sum_{i=2}^na_i^2<\frac12.
\end{equation}
Since $\phi(t)\le t$ for $0\le t<1$, we obtain
\begin{equation}
\sum_{i=1}^n\phi(a_i^2)\le-r\log r+\sum_{i=2}^n a_i^2=r(1-\log r).
\end{equation}
The function $r\mapsto r(1-\log r)$ is increasing on $[0,\frac12]$, and therefore
\begin{equation}
\sum_{i=1}^n\phi(a_i^2)\le\tfrac12(1+\log2).
\end{equation}
Using Lemma~\ref{lem:large-coordinate-1}, we obtain
\begin{equation}
A_{n,p}(a)\le\exp\left\{\big(\tfrac14-\tfrac{1}{2p}\big)(1+\log2)\right\}\le G_p,
\end{equation}
thus completing the proof of \eqref{eq:final-goal-2} for $2<p<4$.
\end{proof}

\medskip


\noindent {\bf Case 2. Suppose that ${\bf 4\leq p\leq26}$.} We shall use the following numerical lemma. 

\begin{lemma}\label{lem:large-coordinate-2}
For every $4\leq p \leq 26$ and every $\frac1{\sqrt2} < a<1$, we have
\begin{equation}\label{eq:BBprime-scalar}
\sup_{q\ge2}\left\{\sqrt{2(1-a^2)}\,q^{1/p}(aG_p)^{q-1}\right\}\ge1.
\end{equation}
\end{lemma}

The proof is deferred to Section \ref{sec:technical}.

\begin{proof} [Proof of Case 2.]
Unless $a_1=1$, in which case the conclusion holds trivially, denote by
\begin{equation} \label{eq:sigma-b-0}
\sigma = \sqrt{1-a_1^2} \qquad \mbox{and} \qquad b_i = \frac{a_i}{\sigma}, \ \mbox{where }  i\in\{2,\ldots,n\}.
\end{equation}
so that $\sum_{i=2}^n b_i^2=1.$ Let $X_1,\ldots,X_n$ be i.i.d.~random variables with density $f_p$ and set $B=\sum_{i=2}^n b_iX_i.$ If $\theta$ denotes the density of $(B+B')/\sqrt2$ and $q\geq2$, then \eqref{eq:holder-2} gives
\begin{equation} \label{eq:use-Holder-426}
A_{n,p}(a)  \leq  2^{-1/2q}a_1^{-(1-1/q)} \sigma^{-1/q} q^{-1/pq} \theta(0)^{1/q},
\end{equation}
Using \eqref{eq:dens-at-0}, we have
\begin{equation} \label{eq:bound-2-426}
\theta(0) = A_{2n-2,p}\big(\tfrac{1}{\sqrt{2}}b,\tfrac{1}{\sqrt{2}}b\big) \leq G_p,
\end{equation}
where the last inequality follows because $\big\|\big(\tfrac{1}{\sqrt{2}}b,\tfrac{1}{\sqrt{2}}b\big)\big\|_\infty \leq \frac{1}{\sqrt{2}}$ and for such vectors we already know inequality \eqref{eq:final-goal-2} by Theorem \ref{thm:Psi} and \eqref{eq:derivation-small-coefficients}. Combining \eqref{eq:use-Holder-426} and \eqref{eq:bound-2-426}, we finally get
\begin{equation}
\bigg(\frac{A_{n,p}(a)}{G_p} \bigg)^q \leq \frac{1}{\sqrt{2}(a_1G_p)^{q-1} \sigma q^{1/p}} = \Big( \sqrt{2(1-a_1^2)}\, q^{1/p}(a_1G_p)^{q-1} \Big)^{-1}.
\end{equation}
The conclusion follows by optimizing over $q\geq2$ and invoking Lemma \ref{lem:large-coordinate-2}.
\end{proof}

\medskip


\noindent {\bf Case 3. Suppose that ${\bf  p>26}$.} We further distinguish several cases depending on the distribution of the coefficient vector $a$.

\smallskip

\noindent {\bf Subcase 3.1. Suppose that $a_1 \geq 2^{\frac{1}{p}-\frac{1}{2}}$.} Then, the last part of Lemma \ref{lem:holder} gives
\begin{equation}
A_{n,p}(a) \leq \frac{1}{a_1} \leq 2^{\frac12-\frac1p}=T_p,
\end{equation}
which implies \eqref{eq:final-goal-2} in this range. \hfill$\Box$

\smallskip

Before distinguishing further cases, denote again by
\begin{equation} \label{eq:sigma-b}
\sigma = \sqrt{1-a_1^2} \qquad \mbox{and} \qquad b_j = \frac{a_j}{\sigma}, \ \mbox{where }  j\in\{2,\ldots,n\}.
\end{equation}
Moreover, let
\begin{equation}
S(b) = \sum_{j=2}^n b_j^4.
\end{equation}
and observe that since $|b_j|\leq 1$ for every $j$, we have $S(b) \leq \sum_{j=2}^n b_j^2 = 1$. 

The proof varies depending on the size of $S(b)$. If $S(b)$ is close to 0, we use the probabilistic representation \eqref{eq:break-exp}, replace the random vector $B$ by a Gaussian, and control the error term by imitating Lindeberg's proof of the central limit theorem. If $S(b)$ is of intermediate size, we repeat the proof of Case 2, except that we now use the deficit term of Theorem \ref{thm:Psi-deficit}. Finally, if $S(b)$ is close to 1, we follow the inductive method of \cite{ENT24},\mbox{ albeit with several technical modifications.}

\smallskip

\noindent {\bf Subcase 3.2. Suppose that $\frac{1}{\sqrt{2}} < a_1 < 2^{\frac{1}{p}-\frac{1}{2}}$ and $S(b) \leq \frac{1}{100}$.} We start with two lemmas.

\begin{lemma} \label{lem:gamma-exp}
Fix $p\geq2$. Then, for every $s>0$,
\begin{equation} \label{eq:gamma-exp}
\big| \gamma_p(s) - e^{-c_p s^2} \big| \leq \tfrac{c_p^2}{2} s^4.
\end{equation}
\end{lemma}

\begin{proof}
The one-sided bound follows readily from
\begin{equation}
\gamma_p(s) \stackrel{\eqref{eq:lb-gammap}}{\geq}1-c_p s^2  \qquad \mbox{and} \qquad e^{-c_p s^2} \leq 1-c_p s^2 + \frac{c_p^2 s^4}{2},
\end{equation}
where the last inequality is a consequence of $e^{-y} \leq 1-y + y^2/2$ for $y\geq0$. For the opposite bound, use the inequality $e^{-c_ps^2} \geq 1-c_ps^2$ along with $\cos t \leq 1-t^2/2+t^4/24$ which gives
\begin{equation}
\gamma_p(s) = \frac{1}{\Gamma(1+\frac1p)} \int_0^\infty \cos(sx) e^{-x^p}\,\diff x \leq 1 - c_ps^2 +\frac{\Gamma(1+\frac5p)}{120\Gamma(1+\frac1p)} s^4.
\end{equation}
To obtain the conclusion, now notice that
\begin{equation}
\frac{\Gamma(1+\frac5p)}{120\Gamma(1+\frac1p)} \stackrel{\eqref{eq:thetap}}{=} (1-\theta_p) \frac{\Gamma(1+\frac3p)^2}{72\Gamma(1+\frac1p)^2} \leq  \frac{\Gamma(1+\frac3p)^2}{72\Gamma(1+\frac1p)^2}
\end{equation}
since, by Lemma \ref{lem:thetap}, we have $\theta_p\geq\theta_2=0$.
\end{proof} 

\begin{lemma} \label{lem:weighted-gaussian-comparison}
Fix $p>26$ and let $\psi:[0,\infty)\to\R$ be either the function $\psi =\gamma_p$ or $\psi(s) = e^{-c_p s^2}$, $s\geq0$, where $c_p$ is given by \eqref{eq:cp}. Then, for every $u\geq8$, we have
\begin{equation} \label{eq:weighted-gaussian-comparison}
\int_0^\infty s^4|\psi(s)|^u\, \diff s \le\tfrac65\int_0^\infty s^4e^{-uc_ps^2}\, \diff s=\frac{9\sqrt\pi}{20c_p^{5/2}}\,u^{-5/2}.
\end{equation}
\end{lemma}

The proof of this lemma is deferred to Section \ref{sec:technical}.

\begin{proof} [Proof of Subcase 3.2]
In this range, $b_j^2\le \sqrt {S(b)}\leq \frac1{10}$ for every $j\in\{2,\ldots,n\}$. Recall that $f_p$ is the probability density $e^{-\beta_p^p |x|^p}$ and let $Z$ be a centered Gaussian random vector with
\begin{equation} \label{eq:mp}
\mathrm{Var}Z = \int_\R x^2 f_p(x) \, \diff x = \frac{\Gamma(1+\frac{3}{p})}{12\Gamma(1+\frac1p)^3} \stackrel{\eqref{eq:cp}}{=} \frac{2c_p}{\beta_p^2}  \eqdef m_p^2.
\end{equation} 
The main idea in this range is to replace the weighted sum in \eqref{eq:break-exp} by the Gaussian variable $Z$. Formula \eqref{eq:section-fourier} combined with a change of variables yields
\begin{equation}
A_{n,p}(a) = \frac{\beta_p}{\pi} \int_0^\infty \prod_{i=1}^n \gamma_p(a_i s)\,\diff s = \frac{\beta_p}{\pi \sigma} \int_0^\infty \gamma_p\big(\tfrac{a_1 s}{\sigma}\big) \prod_{i=2}^n \gamma_p(b_i s)\,\diff s.
\end{equation}
Similarly, by a change of variables and Plancherel's theorem, we get
\begin{equation}
\begin{split}
\E\left[\frac{1}{a_1}f_p\left(\frac{\sigma Z}{a_1}\right)\right] = \frac{1}{a_1}\sqrt{\frac{2}{\pi m_p^2}} \int_0^\infty e^{-\frac{\beta_p^p \sigma^p x^p}{a_1^p}} e^{-\frac{x^2}{2m_p^2}} \, \diff x & = \frac{1}{\sigma \sqrt{\pi c_p}} \int_0^\infty e^{-y^p} e^{-\frac{a_1^2 y^2}{4 c_p \sigma^2}}\, \diff y
\\ & = \frac{\beta_p}{\pi \sigma} \int_0^\infty \gamma_p \big(\tfrac{a_1 s}{\sigma}\big) e^{-c_p s^2}\,\diff s.
\end{split}
\end{equation}
Since $e^{-c_ps^2}= \prod_{i=2}^ne^{-c_pb_j^2s^2}$, we thus get
\begin{equation}	\label{eq:diffuse-fourier-identity}
A_{n,p}(a) - \E\left[\frac{1}{a_1}f_p\left(\frac{\sigma Z}{a_1}\right)\right]  = \frac{\beta_p}{\pi \sigma} \int_0^\infty \gamma_p\big(\tfrac{a_1 s}{\sigma}\big) \bigg( \prod_{i=2}^n \gamma_p(b_is) - \prod_{i=2}^n e^{-c_p b_i^2 s^2}\bigg) \, \diff s.
\end{equation}
Therefore, by telescoping,
\begin{equation}	\label{eq:telescoping}
\begin{split}
\bigg| A_{n,p}(a) -&  \E\left[\frac{1}{a_1}f_p\left(\frac{\sigma Z}{a_1}\right)\right]\bigg|  \leq \frac{\beta_p}{\pi \sigma} \int_0^\infty \big|\gamma_p\big(\tfrac{a_1 s}{\sigma}\big)\big| \bigg| \prod_{i=2}^n \gamma_p(b_is) - \prod_{i=2}^n e^{-c_p b_i^2 s^2}\bigg| \, \diff s
\\ & \leq \frac{\beta_p}{\pi \sigma} \sum_{i=2}^n \int_0^\infty \big|\gamma_p\big(\tfrac{a_1 s}{\sigma}\big)\big| \Big| \gamma_p(b_is) - e^{-c_p b_i^2 s^2} \Big| \prod_{k\neq i} \big| \psi_{i,k}(b_ks)\big| \, \diff s
\end{split}
\end{equation}
where $\psi_{i,k} = \gamma_p$ or $\psi_{i,k}(s) = e^{-c_p s^2}$, $s>0$. Applying inequality \eqref{eq:gamma-exp} of Lemma \ref{lem:gamma-exp}, we thus get
\begin{equation}	\label{eq:telescoping-2}
\begin{split}
\bigg| A_{n,p}(a) -  \E\left[\frac{1}{a_1}f_p\left(\frac{\sigma Z}{a_1}\right)\right]\bigg| \leq \frac{\beta_p c_p^2}{2\pi \sigma} \sum_{i=2}^n b_i^4 \int_0^\infty s^4 \big|\gamma_p\big(\tfrac{a_1 s}{\sigma}\big)\big| \prod_{k\neq i} \big| \psi_{i,k}(b_ks)\big|  \, \diff s
\end{split}
\end{equation}
Hölder's inequality for the measure $s^4\,\diff s$ with weights $b_k^2$ then gives
\begin{equation*}
\begin{aligned}
\int_0^\infty s^4\big|\gamma_p\big(\tfrac{a_1s}{\sigma}\big)\big| \prod_{k\ne i}|\psi_{i,k}(b_ks)|\, \diff s
\le\left(\int_0^\infty s^4\big|\gamma_p\big(\tfrac{a_1s}{\sigma}\big)\big|^{1/b_i^2}\,\diff s\right)^{b_i^2}\prod_{k\ne i}\left(\int_0^\infty s^4\big|\psi_{i,k}(b_ks)\big|^{1/b_k^2}\, \diff s\right)^{b_k^2}.
\end{aligned}
\end{equation*}
Since \(b_k^2<1/10\), all exponents $1/b_k^2$ are at least \(8\), so
Lemma~\ref{lem:weighted-gaussian-comparison} applies.  For \(k\ne i\),
a change of variables gives
\begin{equation}
\int_0^\infty s^4\big|\psi_{i,k}(b_ks)\big|^{1/b_k^2}\, \diff s\le\frac{9\sqrt\pi}{20c_p^{5/2}},
\end{equation}
while
\begin{equation}
\int_0^\infty s^4\big|\gamma_p\big(\tfrac{a_1s}{\sigma}\big)\big|^{1/b_i^2}\,\diff s \le \frac{9\sqrt\pi}{20c_p^{5/2}}\Big(\frac{\sigma b_i}{a_1}\Big)^5
\end{equation}
Consequently,
\begin{equation}
\int_0^\infty s^4\big|\gamma_p\big(\tfrac{a_1s}{\sigma}\big)\big| \prod_{k\ne i}|\psi_{i,k}(b_ks)|\, \diff s
\le \frac{9\sqrt\pi}{20c_p^{5/2}} \Big(\frac{\sigma b_i}{a_1}\Big)^{5b_i^2} \le \frac{9\sqrt\pi}{20c_p^{5/2}},
\end{equation}
since $\sigma \leq a_1$ and $b_i \leq 1$. Substituting into \eqref{eq:telescoping-2} and summing over
\(i\), we get
\begin{equation}
\bigg| A_{n,p}(a) -  \E\left[	\frac{1}{a_1} f_p\left(\frac{\sigma Z}{a_1}\right) \right]\bigg| \le \frac{9\beta_p}{40\sqrt\pi\,\sigma\sqrt{c_p}}S(b).
\end{equation}
Now notice that
\begin{equation} \label{eq:87}
\frac{\beta_p}{\sqrt{c_p}} = \sqrt{\frac{24 \Gamma(1+\frac1p)^3}{\Gamma(1+\frac3p)}} \le\sqrt{24},
\end{equation}
since the logarithm of the left-hand side is $\frac12\log24+\frac32\log\Gamma(1+\frac1p)-\frac12\log\Gamma(1+\frac3p)$, whose derivative with respect to $\frac1p$ equals $\frac32\big((\log\Gamma)'(1+\frac1p)-(\log\Gamma)'(1+\frac3p)\big)\leq0$ because $(\log\Gamma)'$ is increasing, and the limit of the left-hand side as $p\to\infty$ is $\sqrt{24}$. Therefore, we conclude that
\begin{equation} \label{eq:gaussian-comparison}
\bigg| A_{n,p}(a) -  \E\left[ \frac{1}{a_1} f_p\left(\frac{\sigma Z}{a_1}\right) \right]\bigg| \le \frac{9\sqrt{6}}{5\sqrt{7\pi}} S(b) < S(b) \leq \frac{1}{100}.
\end{equation}
We will estimate the Gaussian term.  Set
\begin{equation}
\forall \ x\in\R,\qquad H(x)=\frac{1}{a_1} f_p\left(\frac{\sigma x}{a_1}\right)
\end{equation}
and note that for every $x\in\R$,
\begin{equation}
0\le H(x)\le\frac{1}{a_1}, \qquad \int_{\mathbb R}H(x)\,\diff x=\frac1\sigma.
\end{equation}
Moreover, for 
$R=\frac{a_1}{2\sigma}$, set 
\begin{equation}
H_0(x)=\frac{1}{a_1} \mathbf1_{[-R,R]}(x)
\end{equation}
and observe that for every $x\in\R$,
\begin{equation}
\int_{\mathbb R}H_0(x)\, \diff x = \frac{2R}{a_1} = \frac1\sigma = \int_{\mathbb R}H(x)\, \diff x.
\end{equation}
We apply the bathtub principle \cite[Theorem~1.14]{LL01}: since the density $g$ of the Gaussian random variable $Z$ is even and decreasing on $[0,\infty)$,
\begin{equation}
(H_0(x)-H(x))(g(x)-g(R))\ge0
\end{equation}
for every $x\in\R$ and thus
\begin{equation}
\begin{aligned}
	\int_{\mathbb R}(H_0-H)g  &= \int_{\mathbb R}(H_0-H)(g-g(R)) + g(R)\int_{\mathbb R}(H_0-H)\\
	&= \int_{\mathbb R}(H_0-H)(g-g(R)) \ge0.
\end{aligned}
\end{equation}
Hence
\begin{equation}
\begin{aligned}
\E\left[\frac{1}{a_1} f_p\left(\frac{\sigma Z}{a_1}\right)\right]=\int_{\mathbb R}H(x)g(x)\,\diff x &  \le\int_{\mathbb R}H_0(x)g(x)\,\diff x
 \\ & =\frac1a_1 \mathbb P\left( |Z|\le\frac{a_1}{2\sigma}\right)=\frac{1}{a_1}\left(2\Phi\left(\frac{a_1}{2\sigma m_{p}}\right)-1\right),
\end{aligned}
\end{equation}
where $\Phi$ is the Gaussian distribution function and $m_p$ is given by \eqref{eq:mp}. Denoting by $r_0=a_1/\sigma \geq1$ and using \(m_{p}\ge1/\sqrt{12}\) which follows from \eqref{eq:87}, we can further bound this as
\begin{equation}
\E\left[ \frac{1}{a_1} f_p\left(\frac{\sigma Z}{a_1}\right)\right] \leq \frac{1}{a_1}\left(2\Phi\left(\frac{\sqrt{3}a_1}{\sigma}\right)-1\right) = \frac{\sqrt{1+r_0^2}}{r_0}\bigl(2\Phi(\sqrt3\,r_0)-1\bigr).
\end{equation}
This function is decreasing for \(r_0\ge1\).  Indeed, its derivative is nonpositive precisely when
\begin{equation}
2\Phi(\sqrt3\,r)-1 \ge 2\sqrt3\,r(1+r^2)\Phi'(\sqrt3\,r),
\end{equation}
and the difference between the two sides has derivative
\begin{equation}
6\sqrt3\,r^4\Phi'(\sqrt3\,r)>0
\end{equation}
and is positive at $r=1$.  Hence
\begin{equation}
\E\left[ \frac1{a_1} f_p\left(\frac{\sigma Z}{a_1}\right) \right] \le \sqrt2\,\big(2\Phi(\sqrt3)-1\big) < \frac43.
\end{equation}
Combining this with \eqref{eq:gaussian-comparison},
\begin{equation}
A_{n,p}(a) < \frac43+\frac1{100} < \frac{11}{8} < 2^{\frac12-\frac{1}{26}} \leq T_p
\end{equation}
for $p\geq26$. This concludes the proof of \eqref{eq:final-goal-2} for this case.
\end{proof}

\smallskip

\noindent {\bf Subcase 3.3. Suppose that $\frac{1}{\sqrt{2}} < a_1 < 2^{\frac{1}{p}-\frac{1}{2}}$, $\frac{1}{100}<S(b)$ and $S(b)(1-\sqrt{S(b)}) \geq \tfrac{50}{9p^2}$.} We start with a numerical lemma.

\begin{lemma} \label{lem:3.3}
For every $p>26$ and every $\frac{1}{\sqrt2}< a < 1$, we have
\begin{equation} \label{eq:lem-3.3}
\sup_{q\geq2} \Big\{ \sqrt{2(1-a^2)} \big( a \max\{T_p,G_p\}\big)^{q-1} q^{1/p} \Big\} \geq \exp\Big(-\tfrac{1}{12p^2}\Big).
\end{equation}
\end{lemma}

The proof is deferred to Section \ref{sec:technical}.

\begin{proof} [Proof of Subcase 3.3]
We proceed with the H\"older argument, as in the range $4\leq p\leq26$. Combining \eqref{eq:holder-2} and \eqref{eq:bound-2-426}, for every $q\geq2$, we get
\begin{equation}
A_{n,p}(a) \leq 2^{-1/2q} a_1^{-(1-1/q)} \sigma^{-1/q} q^{-1/pq} A_{2n-2,p}\big(\tfrac{1}{\sqrt2}b,\tfrac{1}{\sqrt2}b\big)^{1/q}.
\end{equation}
For the latter section function, we use \eqref{eq:holder-Psi} to bound it by
\begin{equation}
A_{2n-2,p}\big(\tfrac{1}{\sqrt2}b,\tfrac{1}{\sqrt2}b\big) \leq  \tfrac{2}{\pi} \Gamma\big(1+\tfrac1p\big) \prod_{i=2}^n \Psi_p\big( 2/b_i^2\big)^{b_i^2}.
\end{equation}
Using the last identity of \eqref{eq:derivation-small-coefficients}, we can rewrite this as
\begin{equation} \label{eq:used-deficit-largep}
\begin{split}
\frac{A_{2n-2,p}\big(\tfrac{1}{\sqrt2}b,\tfrac{1}{\sqrt2}b\big)}{\max\{T_p,G_p\}} & \leq \prod_{i=2}^n \bigg( \frac{\Psi_p\big( 2/b_i^2\big)}{\max\{\Psi_p(2),\Psi_p(\infty)\}} \bigg)^{b_i^2} \stackrel{\eqref{eq:Psi-bound-deficit}}{\leq} \prod_{i=2}^n \Big( 1 - \tfrac{3}{200} (b_i^2-b_i^4)  \Big)^{b_i^2} 
\\ & \leq \exp\bigg( - \frac{3}{200} \sum_{i=2}^n b_i^4(1-b_i^2)\bigg)
\end{split}
\end{equation}
where the second inequality used the sharp upper bound for $\Psi_p$ with a deficit proven in Theorem \ref{thm:Psi-deficit}. By monotonicity of the $\ell_p$ norms,
\begin{equation}
\sum_{i=2}^n b_i^6 \leq \Big(\sum_{i=2}^n b_i^4\Big)^{3/2} = S(b)^{3/2}
\end{equation}
and thus \eqref{eq:used-deficit-largep} gives
\begin{equation}
\frac{A_{2n-2,p}\big(\tfrac{1}{\sqrt2}b,\tfrac{1}{\sqrt2}b\big)}{\max\{T_p,G_p\}}  \leq \exp\Big( - \tfrac{3}{200}  S(b)\big(1-\sqrt{S(b)}\big)\Big) \leq \exp\Big( - \tfrac{1}{12p^2} \Big)
\end{equation}
by the assumption on $S(b)$. Combining all the above, we thus get
\begin{equation}
\left(\frac{A_{n,p}(a)}{\max\{T_p,G_p\}}\right)^q \leq \exp\Big( - \tfrac{1}{12p^2} \Big) \Big( \sqrt{2(1-a_1^2)} \big( a_1 \max\{T_p,G_p\}\big)^{q-1} q^{1/p} \Big)^{-1}.
\end{equation}
The conclusion follows by optimizing over $q\geq2$ and invoking Lemma \ref{lem:3.3}.
\end{proof}

\smallskip

\noindent {\bf Subcase 3.4. Suppose that $\frac{1}{\sqrt{2}} < a_1 < 2^{\frac{1}{p}-\frac{1}{2}}$, $\frac{49}{50} < S(b)$ and $S(b)(1-\sqrt{S(b)})<\frac{50}{9p^2}$.} Before passing to the proof, let us point out that for $p\geq p_0$ this is the conceptually more delicate case, as for the extremizer $a=\frac{e_1+e_2}{\sqrt{2}}$, we have $b=(1,0,\ldots,0)$ and thus $S(b)=1$. Moreover, observe that this is indeed the only remaining case as assuming that $S(b)(1-\sqrt{S(b)})<\frac{50}{9p^2}$ together with $S(b)>\tfrac{1}{100}$ already implies $S(b)> \frac{49}{50}$ for $p>26$. We also assume that $a_1^2+a_2^2<1$, as otherwise we are reduced to the planar case $n=2$.

\begin{proof} [Proof of Subcase 3.4] We shall use induction on the dimension $n$, along the lines of \cite{ENT24}. The cases $n\in\{1,2\}$ are clear by direct computation. Following the notation of Lemma \ref{lem:min-bound}, let $R_1,\ldots,R_n$ be i.i.d.~random variables with density $\alpha_p^{-1} x^p e^{-x^p} {\bf 1}_{x>0}$, where $\alpha_p = \frac{1}{p}\Gamma(1+\frac{1}{p})$, and $\xi_1,\ldots,\xi_n$ be i.i.d.~random vectors uniformly distributed on the unit sphere $\mb{S}^2$, independent of the random variables $R_i$. If $X=a_1R_1\xi_1+a_2R_2\xi_2$ and $Y=\sum_{i\geq3} a_i R_i \xi_i$, then Lemma \ref{lem:min-bound} gives
\begin{equation}
A_{n,p}(a) \leq \Gamma\big(1+\tfrac1p\big) \mb{E} \big[ \min \big\{ |X|^{-1}, \mb{E}|Y|^{-1}\big\} \big].
\end{equation}
By the inductive hypothesis,
\begin{equation}
\mb{E}|Y|^{-1} = \frac{1}{\sqrt{1-a_1^2-a_2^2}} \mb{E} \Bigg| \frac{\sum_{i\geq3} a_i R_i \xi_i}{\sqrt{1-a_1^2-a_2^2}} \Bigg|^{-1} \stackrel{\eqref{eq:spherical}}{\leq} \frac{\max\{T_p,G_p\}}{\Gamma(1+\frac1p) \sqrt{1-a_1^2-a_2^2}}.
\end{equation}
Let
\begin{equation}
\beta \eqdef
\frac{\Gamma(1+\frac1p)}{\max\{T_p,G_p\}} \sqrt{1-a_1^2-a_2^2} \quad \mbox{and} \quad \alpha\eqdef
\frac{\Gamma(1+\frac1p)}{T_p} \sqrt{1-a_1^2-a_2^2}.
\end{equation}
Therefore, since $\beta\leq \alpha$  and $\min\{\chi,\psi\} = \chi- (\chi-\psi)_+$,
\begin{equation} \label{eq:34-break}
\begin{split}
A_{n,p}(a) & \leq \Gamma\big(1+\tfrac1p\big) \mb{E}\min\big\{|X|^{-1},\beta^{-1}\big\} 
\\ & \leq  \Gamma\big(1+\tfrac1p\big) \mb{E}\min\big\{|X|^{-1},\alpha^{-1}\big\} + \Gamma\big(1+\tfrac1p\big) (\beta^{-1}-\alpha^{-1}) \mb{P}\{|X|\leq \alpha\}
\\ & = \Gamma\big(1+\tfrac1p\big) \mb{E} \big[ |X|^{-1} - \big(|X|^{-1}-\alpha^{-1}\big)_+ \big] + \frac{\max\{T_p,G_p\}-T_p}{\sqrt{1-a_1^2-a_2^2}} \ \mb{P}\{|X|\leq \alpha\}
\\ & = \frac{1}{\|(a_1,a_2)\|_p} - \Gamma\big(1+\tfrac1p\big) \mb{E} \big(|X|^{-1}-\alpha^{-1}\big)_+ \frac{\max\{T_p,G_p\}-T_p}{\sqrt{1-a_1^2-a_2^2}} \ \mb{P}\{|X|\leq \alpha\},
\end{split}
\end{equation}
where the computation of $\mb{E}|X|^{-1}$ in the last equality follows from \eqref{eq:spherical}. 

We start by treating the latter term, which is nonzero only for $26<p<p_0$. For fixed $a,b>0$,
\begin{equation}
|a\xi_1+b\xi_2|^2 = a^2+b^2+2ab\langle\xi_1,\xi_2\rangle.
\end{equation}
Since $\langle\xi_1,\xi_2\rangle$ is uniform on $[-1,1]$, for $t\geq0$, we get
\begin{equation}
\P_{\xi_1,\xi_2}\big\{|a\xi_1+b\xi_2|<t\big\} =
\begin{cases}
0, & t\le |a-b|,\\[1mm]
\dfrac{t^2-(a-b)^2}{4ab},
& |a-b|<t<a+b,\\[2mm]
1, & t\ge a+b,
\end{cases}
\end{equation}
or equivalently,
\begin{equation}
\P_{\xi_1,\xi_2}\big\{|a\xi_1+b\xi_2|<t\big\} = \frac{\bigl(t^2-(a-b)^2\bigr)_+ - \bigl(t^2-(a+b)^2\bigr)_+}{4ab}.
\end{equation}
Averaging over \(R_1,R_2\), we obtain
\begin{equation} \label{eq:ball-formula}
\P\big\{|X|<t\big\} = \frac1{4a_1a_2} \E\left[ \frac{\bigl(t^2-(a_1R_1-a_2R_2)^2\bigr)_+-\bigl(t^2-(a_1R_1+a_2R_2)^2\bigr)_+}{R_1R_2}\right].
\end{equation}
Recall that \(R_i\) have density $\frac{p}{\Gamma(1+\frac1p)}\,r^p e^{-r^p}\textbf{1}_{r\geq0},$
hence for every measurable function \(F\),
\begin{equation}
\E\left[\frac{F(R_1,R_2)}{R_1R_2}\right] = \frac{p^2}{\Gamma(1+\frac1p)^2} \int_0^\infty\int_0^\infty F(r_1,r_2) (r_1r_2)^{p-1} e^{-r_1^p-r_2^p} \, \diff r_1 \diff r_2.
\end{equation}
The function $r\mapsto p r^{p-1}e^{-r^p} {\bf 1}_{r\geq0}$ integrates to 1 and thus, if \(S_1,S_2\) are independent random variables with this density, then
\begin{equation}\label{eq:R-to-W}
\E\left[\frac{F(R_1,R_2)}{R_1R_2}\right] = \frac{1}{\Gamma(1+\frac1p)^2} \ \E F(S_1,S_2).
\end{equation}
Therefore, combining \eqref{eq:ball-formula} and \eqref{eq:R-to-W} yields
\begin{equation} \label{eq:ball-formula-2}
\P\big\{|X|<t\big\} = \frac{1}{4a_1a_2\,\Gamma(1+\frac1p)^2} \left( \E\bigl(t^2-(a_1S_1-a_2S_2)^2\bigr)_+ - \E\bigl(t^2- (a_1S_1+a_2S_2)^2\bigr)_+ \right).
\end{equation}
In particular, this trivially yields 
\begin{equation} \label{eq:ball-bound}
\frac{\mb{P}\{|X|\leq \alpha\}}{\sqrt{1-a_1^2-a_2^2}} \leq \frac{\alpha^2}{4a_1a_2\sqrt{1-a_1^2-a_2^2}\Gamma(1+\frac1p)^2} = \frac{\sqrt{1-a_1^2-a_2^2}}{4a_1a_2 T_p^2}.
\end{equation}
To estimate the remaining numerical constant observe that
\begin{equation} \label{eq:largee-b}
\frac{49}{50} \leq S(b) =\sum_{j=2}^nb_j^4 \leq b_2^2 \sum_{j=2}^n b_j^2 = b_2^2,
\end{equation}
i.e. $a_2^2 \geq \frac{49}{50} (1-a_1^2)$. Since $a_1> \frac{1}{\sqrt2}$, this gives
\begin{equation}
\frac{\sqrt{1-a_1^2-a_2^2}}{4a_1a_2} \leq \frac{\sqrt{\frac{1}{50}(1-a_1^2)}}{4a_1 \sqrt{\frac{49}{50} (1-a_1^2)}} < \frac{1}{4\cdot \frac{1}{\sqrt2} \cdot \sqrt{49}} <1 
\end{equation}
Combining this with $T_p \geq1$ and \eqref{eq:ball-bound} and plugging into \eqref{eq:34-break}, we thus get
\begin{equation} \label{eq:34-almost-done}
A_{n,p}(a) \leq \frac{1}{\|(a_1,a_2)\|_p} - \Gamma\big(1+\tfrac1p\big) \mb{E} \big(|X|^{-1}-\alpha^{-1}\big)_+ + \big(\max\{T_p,G_p\}-T_p\big).
\end{equation}

We are now left with estimating the middle term in this expression. This is the content of the next proposition, which was previously proven in \cite[Proposition~15]{ENT24} for $p\geq10^{15}$.

\begin{proposition} \label{prop:ent}
Under the notation and assumptions above, for $p>26$, we have
\begin{equation} \label{eq:prop-ent}
 \mb{E} \big(|X|^{-1}-\alpha^{-1}\big)_+ \geq \frac{1}{\Gamma(1+\frac1p)} \bigg( \frac{1}{\|(a_1,a_2)\|_p} - T_p \bigg).
\end{equation}
\end{proposition}

\begin{proof}
If the right-hand side is nonpositive, there is nothing to prove. Otherwise,
\begin{equation} \label{eq:p<T}
\|(a_1,a_2)\|_p < T_p^{-1}.
\end{equation}
By comparison of the \(\ell_p\) and \(\ell_2\) norms in $\R^2$,
$\|(a_1,a_2)\|_p \geq T_p^{-1} \|(a_1,a_2)\|_2,$
and thus
\begin{equation*}
\begin{split}
\frac{1}{\Gamma(1+\frac1p)} \Big( \frac{1}{\|(a_1,a_2)\|_p} - T_p\Big) & \leq 
\frac{T_p}{\Gamma(1+\frac1p)} \bigg( \frac{1}{\sqrt{a_1^2+a_2^2}}-1 \bigg) = \frac{T_p}{\Gamma(1+\frac1p)}\frac{1-a_1^2-a_2^2}{\sqrt{a_1^2+a_2^2}(1+\sqrt{a_1^2+a_2^2})}
\\ & = \bigg( \frac{T_p}{\Gamma(1+\frac1p)} \bigg)^3 \frac{\alpha^2}{\sqrt{a_1^2+a_2^2}(1+\sqrt{a_1^2+a_2^2})}.
\end{split}
\end{equation*}
Combining the inequality $\Gamma(1+\frac1p) \geq 2^{-\frac1p}$ which holds by log-convexity of $\Gamma$, as $\log\Gamma(1+\frac1p) \geq \frac1p (\log\Gamma)'(1) > -\frac{\log2}{p}$, with the bound $a_1^2+a_2^2 \geq \frac{49}{50}$ which follows from \eqref{eq:largee-b}, allows us to estimate  the numerical constant as
\begin{equation}
\bigg( \frac{T_p}{\Gamma(1+\frac1p)} \bigg)^3 \frac{1}{\sqrt{a_1^2+a_2^2}(1+\sqrt{a_1^2+a_2^2})} \leq \sqrt{2}^3 \frac{1}{\sqrt{\frac{49}{50}}(1+\sqrt{\frac{49}{50}})} < \frac32.
\end{equation}
Therefore, the proof shall be complete once we prove that
\begin{equation} \label{eq:3ov2}
 \mb{E} \big(|X|^{-1}-\alpha^{-1}\big)_+ \geq \tfrac{3}{2} \alpha^2.
\end{equation}

By the layer-cake formula,
\begin{equation}\label{eq:layer-cake}
\E\bigl(|X|^{-1}-\alpha^{-1}\bigr)_+ =
\int_0^\alpha\frac{\P\{|X|<t\}}{t^2}\, \diff t.
\end{equation}
Let $S_1, S_2$ be independent random variables with density $r\mapsto p r^{p-1} e^{-r^p} {\bf 1}_{r\geq0}$ and denote by $\Delta =a_1S_1-a_2S_2$ and $\Sigma = a_1S_1+a_2S_2$. Then, \eqref{eq:ball-formula-2} gives
\begin{equation} \label{eq:D-S}
\E\bigl(|X|^{-1}-\alpha^{-1}\bigr)_+ = \frac{1}{4a_1a_2\Gamma(1+\frac1p)^2} 
\int_0^\alpha \E\big( t^2-\Delta^2\big)_+ - \E\big( t^2 - \Sigma^2\big)_+ \, \frac{\diff t}{t^2}
\end{equation}
Since each $S_i$ is log-concave, it is a standard consequence of the Pr\'ekopa--Leindler inequality (see, e.g., \cite[Corollary~1.4.2]{AGM15}) that the same holds for $\Delta$ and $\Sigma$. 

To estimate the first term in \eqref{eq:D-S}, we repeat the proof of \cite[Lemma~12.2.6]{BGVV14}. Let $q:\R\to\R_+$ be the density of the random variable $\Delta$. A direct computation gives
\begin{equation}\label{eq:q-zero}
q(0) = \frac{p^2}{(a_1 a_2)^{p}} \int_0^\infty r^{2p-2} e^{-(\frac{1}{a_1^p}+\frac{1}{a_2^p}) r^p} \, \diff r = p \Gamma\big(2-\tfrac1p\big) \frac{(a_1a_2)^{p-1}}{(a_1^p+a_2^p)^{2-1/p}}.
\end{equation}
Moreover, since \(S_1^p,S_2^p\) are independent standard exponential random
variables,
\begin{equation}\label{eq:q-side-masses}
M_+ \eqdef \int_0^\infty q(s)\,\diff s = \mb{P}\big\{ S_1^p \geq \tfrac{a_2^p}{a_1^p} S_2^p\big\} = \frac{a_1^p}{a_1^p+a_2^p} \quad \mbox{and} \quad M_-\eqdef \int_{-\infty}^0q(s)\,\diff s  = \frac{a_2^p}{a_1^p+a_2^p}.
\end{equation}
Again by the Pr\'ekopa--Leindler inequality, the tail $s\longmapsto\int_s^\infty q(u)\, \diff u$
is log-concave on \([0,\infty)\). Its value and derivative at the origin are \(M_+\) and \(-q(0)\), respectively, hence concavity of the logarithm gives for every $s>0$,
\begin{equation}
\int_s^\infty q(u)\,\diff u \le M_+\exp\left(-\frac{q(0)}{M_+}s\right).
\end{equation}
Applying the same argument to \(q(-u)\), we obtain
\begin{equation}\label{eq:D-small-ball}
\P\big\{|\Delta|<s\big\} \ge M_+\left(1-e^{-sq(0)/M_+}\right) + M_-\left(1-e^{-sq(0)/M_-}\right).
\end{equation}
To estimate the first term in \eqref{eq:D-S}, write
\begin{equation}
\begin{aligned}
\E\big(t^2-\Delta^2\big)_+ &= \int_0^{t^2} \P\big\{ |\Delta|<\sqrt{t^2-u}\big\}\, \diff u
\\ & \stackrel{\eqref{eq:D-small-ball}}{\ge} M_+\int_0^{t^2} \left( 1-e^{-\frac{q(0)}{M_+}\sqrt{t^2-u}} \right)\diff u +
M_-\int_0^{t^2} \left( 1-e^{-\frac{q(0)}{M_-}\sqrt{t^2-u}} \right)\,\diff u.
\end{aligned}
\end{equation}
After the change of variables \(u=t^2(1-s^2)\), this becomes
\begin{equation}
\begin{aligned}
\E\big(t^2-\Delta^2\big)_+ \ge 2t^2\Bigg[ M_+\int_0^1 s\left(1-e^{-tsq(0)/M_+}\right)\, \diff s + M_-\int_0^1 s\left(1-e^{-tsq(0)/M_-}\right)\, \diff s \Bigg].
\end{aligned}
\end{equation}
Integrating with respect to \(t\), we obtain
\begin{equation}
\begin{aligned}
\int_0^\alpha \E\big( t^2-\Delta^2\big)_+ \,\frac{\diff t}{t^2} \ge 2\alpha\Bigg[ M_+\int_0^1 s& \left( 1- \frac{1-\exp\!\left(-\frac{\alpha q(0)}{M_+}s\right)} {\frac{\alpha q(0)}{M_+}s} \right)\,\diff s
\\ &+ M_-\int_0^1s\left(1- \frac{1-\exp\!\left(-\frac{\alpha q(0)}{M_-}s\right)} {\frac{\alpha q(0)}{M_-}s} \right)\,\diff s\Bigg].
\end{aligned}
\end{equation}
Consider the function $h:(0,\infty)\to\R$ given by
\begin{equation}
\forall \ z>0,\qquad h(z) \eqdef \frac{2}{z} \int_0^1s\left(1-\frac{1-e^{-zs}}{zs}\right)\, \diff s = 2\int_0^1 s^2\int_0^1 (1-u) e^{-zsu}\, \diff u \, \diff s,
\end{equation}
which is evidently decreasing and convex on \((0,\infty)\). Therefore,
\begin{equation}
\begin{split}\label{eq:D-exact-integrated}
\int_0^\alpha \E\big( t^2-\Delta^2\big)_+ \,\frac{\diff t}{t^2}  & \ge q(0)\alpha^2 \left[ h\left(\frac{\alpha q(0)}{M_+} \right) + h\left(\frac{\alpha q(0)}{M_-}\right) \right]
\\ &\geq 2q(0)\alpha^2h\left(\frac{\alpha q(0)}{2} \left(\frac{1}{M_+}+\frac{1}{M_-}\right)\right).
\end{split}
\end{equation}
We shall also use the following claim.

\medskip

\noindent {\bf Claim.} Under the assumptions above, we have
\begin{equation}\label{eq:z-sum}
\alpha q(0)\left( \frac{1}{M_+}+\frac{1}{M_-} \right) <2.
\end{equation}

\begin{proof}
From the formulas \eqref{eq:q-zero} and \eqref{eq:q-side-masses},
\begin{equation} \label{eq:alphaq}
\alpha q(0)\left(\frac1{M_+}+\frac1{M_-} \right) = p \sqrt{1-a_1^2-a_2^2}  \frac{\Gamma(1+\frac1p)\Gamma\big(2-\tfrac1p\big)}{T_p}\frac{(a_1^p+a_2^p)^{1/p}}{a_1a_2}.
\end{equation}
Since $S(b)(1-\sqrt{S(b)})<\frac{50}{9p^2}$, we have $1-\sqrt{S(b)}<\frac{50}{9p^2S(b)}$ and therefore
\begin{equation}
1-S(b)=(1-\sqrt{S(b)})(1+\sqrt{S(b)})<\frac{100}{9p^2S(b)}.
\end{equation}
On the other hand, since $\sum_{j\geq3}b_j^4\leq\big(\sum_{j\geq3}b_j^2\big)^2=(1-b_2^2)^2$, we have
\begin{equation}
S(b)\leq b_2^4+(1-b_2^2)^2=1-2b_2^2(1-b_2^2),
\end{equation}
that is, $1-b_2^2\leq\frac{1-S(b)}{2b_2^2}$. Moreover, $S(b)\leq b_2^2\sum_{j\geq2} b_j^2=b_2^2$, thus
\begin{equation}
\frac{1-b_2^2}{2}\leq\frac{1-S(b)}{4S(b)}<\frac{25}{9p^2S(b)^2}.
\end{equation}
Since $S>\frac{49}{50}$, we thus have
\begin{equation}
\frac{1-b_2^2}{2}<\frac{25}{9p^2}\cdot\frac{2500}{2401}<\frac{29}{10p^2},
\end{equation}
which implies that
\begin{equation} \label{eq:bound-b2}
p\sqrt{1-a_1^2-a_2^2} = p\sqrt{1-a_1^2} \sqrt{1-b_2^2} < p\sqrt{\frac{1-b_2^2}{2}} < \sqrt{\frac{29}{10}}.
\end{equation}
Combining this with the fact that $\Gamma(x)<1$ for $x\in(1,2)$ and \eqref{eq:p<T}, \eqref{eq:alphaq} becomes
\begin{equation} \label{eq:bor}
\alpha q(0)\left(\frac1{M_+}+\frac1{M_-} \right) < \sqrt{\frac{29}{10}} \frac{1}{T_p^2a_1a_2}< \sqrt{\frac{29}{10}} \frac{2^{\frac{1}{13}}}{2a_1a_2}.
\end{equation}
Finally, using \eqref{eq:bound-b2} and $a_1\leq 2^{\frac1p-\frac12}$,
\begin{equation} \label{eq:a1a2}
\begin{split}
(2a_1a_2)^2 & = (a_1^2+a_2^2)^2 - (a_1^2-a_2^2)^2 = (a_1^2+a_2^2)^2 - \big(2a_1^2 - (a_1^2+a_2^2)\big)^2 
\\ & > \big(1-\tfrac{29}{10p^2}\big)^2 - \big( 2^{\frac2p} - 1+\tfrac{29}{10p^2}\big)^2 > 1-\frac{21}{2p^2} > \frac{1331}{1352}.
\end{split}
\end{equation}
Plugging this into \eqref{eq:bor} proves the claim.
\end{proof}

\smallskip

\noindent {\it Finishing the proof of Proposition \ref{prop:ent}.} Using the claim and the fact that $h$ is decreasing, we get from \eqref{eq:D-exact-integrated} that
\begin{equation} \label{eq:D-bound}
\begin{split}
\frac{1}{4a_1a_2\Gamma(1+\frac1p)^2} \int_0^\alpha \E\big( t^2-\Delta^2\big)_+ \,\frac{\diff t}{t^2} &\geq \frac{2q(0)\alpha^2h(1)}{4a_1a_2\Gamma(1+\frac1p)^2} = \big(\tfrac12-\tfrac1e\big) \frac{q(0)\alpha^2}{a_1a_2\Gamma(1+\tfrac1p)^2} 
\\ & = \big(\tfrac12-\tfrac1e\big)  \frac{\Gamma(2-\tfrac1p)}{\Gamma(1+\tfrac1p)^2} \frac{(a_1a_2)^{p-2}}{(a_1^p+a_2^p)^{2-1/p}} p\alpha^2
\end{split}
\end{equation}
Bernoulli's inequality and \eqref{eq:a1a2} imply
\begin{equation}\label{eq:a1a2-balance}
\begin{aligned}
(2a_1a_2)^{p-2}> \left(1-\frac{21}{2p^2}\right)^{\frac{p-2}{2}} \ge 1-\frac{21(p-2)}{4p^2} > \frac45,
\end{aligned}
\end{equation}
while \eqref{eq:p<T} gives
\begin{equation} 
(a_1^p+a_2^p)^{2-1/p} = \|(a_1,a_2)\|_p^{2p-1} < 2^{\frac{2p-1}{p} - \frac{2p-1}{2}} = \frac{2^{\frac12-\frac1p}}{2^{p-2}} < \frac{\sqrt2}{2^{p-2}}.
\end{equation}
Finally, the function $\varphi(t) = \log \Gamma(2-t) - 2\log\Gamma(1+t)$, where $t\in(0,\tfrac{1}{26})$, satisfies
\begin{equation}
\varphi''(t) = (\log \Gamma)''(2-t) - 2 (\log\Gamma)''(1+t) \leq (\log\Gamma)''\big(\tfrac{51}{26}\big) - 2(\log\Gamma)''\big(\tfrac{27}{26}\big) < 0 
\end{equation}
by the classical formulas \cite[6.4.10]{AS64} for the derivatives of the digamma function. Therefore, $\varphi$ is concave and $\varphi'(\tfrac{1}{26}) \approx  0.63$ so $\varphi$ is also  increasing. Thus $\varphi(t) \geq \varphi(0) = 0$ which can be equivalently written as $\Gamma(2-\frac1p) > \Gamma(1+\frac1p)^2$ for $p>26$. Plugging these into \eqref{eq:D-bound}, we get
\begin{equation}
\frac{1}{4a_1a_2\Gamma(1+\frac1p)^2} \int_0^\alpha \E\big( t^2-\Delta^2\big)_+ \,\frac{\diff t}{t^2} \geq \big(\tfrac12-\tfrac1e\big) \cdot 1 \cdot \frac{\frac{4}{5}}{\sqrt2} \cdot 26 \cdot \alpha^2 > \tfrac{14}{9} \alpha^2.
\end{equation}

We will now estimate the second term of \eqref{eq:D-S}, showing that
\begin{equation}
\frac{1}{4a_1a_2\Gamma(1+\frac1p)^2} \int_0^\alpha \E\big(t^2-\Sigma^2\big)_+\, \frac{\diff t}{t^2} <\tfrac1{100}\alpha^2.
\end{equation}
For every $u>0$, we have
\begin{equation}
\P\{S_i<u\}= \P\big\{  S_i^p < u^p \big\} = 1-e^{-u^p}\le u^p.
\end{equation}
Therefore, for $t\in(0,\alpha)$, since $\Sigma=a_1S_1+a_2S_2$,
\begin{equation}
\P\{\Sigma<t\} \leq \mb{P}\big\{S_1 < \tfrac{t}{a_1} \mbox{ and } S_2< \tfrac{t}{a_2}\big\} \le \frac{t^{2p}}{(a_1a_2)^p},
\end{equation}
which implies
\begin{equation}
\E\big(t^2-\Sigma^2\big)_+ \le t^2\P\{\Sigma<t\} \le \frac{t^{2p+2}}{(a_1a_2)^p}.
\end{equation}
It follows that
\begin{equation}
\begin{aligned}
\frac{1}{4a_1a_2\Gamma(1+\frac1p)^2} \int_0^\alpha \E\big(t^2-\Sigma^2\big)_+\, \frac{\diff t}{t^2}  &\le\frac{\alpha^{2p+1}}{4(a_1a_2)^{p+1}\Gamma(1+\frac1p)^2(2p+1)}
\\&=\alpha^2\frac{\alpha}{(2a_1a_2)^2\Gamma(1+\frac1p)^2(2p+1)} \left(\frac{\alpha^2}{a_1a_2}\right)^{p-1}.
\end{aligned}
\end{equation}
The bounds \eqref{eq:bound-b2}, \eqref{eq:a1a2} yield $\alpha<\frac{1}{20}$ and $a_1a_2>\frac38$. Combining these with $\Gamma(1+\frac1p)>\frac12$ gives
\begin{equation}
\frac{1}{4a_1a_2\Gamma(1+\frac1p)^2} \int_0^\alpha \E\big(t^2-\Sigma^2\big)_+\, \frac{\diff t}{t^2} < \left(\frac1{400}\cdot \frac83 \right)^{p-1}\alpha^2 < \tfrac1{100}\alpha^2.
\end{equation}
Putting everything together,
\begin{equation}
\begin{aligned}
\E\bigl(|X|^{-1}-\alpha^{-1}\bigr)_+ > \left(\tfrac{14}{9}-\tfrac1{100}\right)\alpha^2 > \tfrac32\alpha^2,
\end{aligned}
\end{equation}
which completes the proof of \eqref{eq:3ov2}, which in turn implies Proposition \ref{prop:ent}.
\end{proof}

\smallskip

\noindent {\it Finishing the proof of Subcase 3.4.} Plugging \eqref{eq:prop-ent} into \eqref{eq:34-almost-done} gives
\begin{equation*}
A_{n,p}(a) \leq \frac{1}{\|(a_1,a_2)\|_p} - \bigg( \frac{1}{\|(a_1,a_2)\|_p} - T_p \bigg)+ \big(\max\{T_p,G_p\}-T_p\big) = \max\{T_p,G_p\}. \hfill\qedhere
\end{equation*}
\end{proof}


\section{Proofs of the remaining technical estimates} \label{sec:technical}

In this section, we prove the lemmas which were omitted from the proof of Theorem \ref{thm:main}.

\begin{proof} [Proof of Lemma \ref{lem:large-coordinate-1}]
Setting $a=\frac1p\in(\frac14,\frac12)$, the desired inequality is equivalent to
\begin{equation}
3\log\Gamma(1+a) - \log\Gamma(1+3a) +\log\tfrac{6}{\pi} - \big(\tfrac12-a\big)(1+\log 2) \geq 0.
\end{equation}
Denoting by $k(a)$ this expression and using the classical formula \cite[6.4.10]{AS64} for the derivatives of the polygamma function, we have
\begin{equation}
k''(a) = 3  \sum_{n=1}^\infty \Big( \frac{1}{(n+a)^2} - \frac{3}{(n+3a)^2} \Big).
\end{equation}
For every $a\leq 1/2$ and $n\in\N$, $(n+3a)^2 < 3(n+a)^2$ and thus $k$ is concave on $[\frac14,\frac12]$. Therefore, $k(a) \geq \min\{k(1/2), k(1/4)\} = 0$, which concludes the proof.
\end{proof}

\begin{proof} [Proof of Lemma \ref{lem:large-coordinate-2}]
If $aG_p\geq1$, the limit of the function $\alpha_p(q) \eqdef \sqrt{2(1-a^2)}q^{1/p}(aG_p)^{q-1}$ at infinity is  $\lim_{q\to\infty} \alpha_p(q)=\infty$ and the claim is obvious. If $aG_p<1$, observe that
\begin{equation}
\alpha_p'(q) = \frac{\sqrt{2(1-a^2)}}{p}q^{1/p-1} (aG_p)^{q-1} \big(1 + pq \log(aG_p) \big).
\end{equation}
which vanishes at $q_p^* = \frac{1}{p\log(1/aG_p)}$. If $q_p^*\leq 2$, that is $aG_p \leq e^{-1/2p}$, then
\begin{equation} \label{eq:61}
\sup_{q\geq2} \alpha_p(q) = \alpha_p(2) = 2^{\frac12+\frac1p}  \sqrt{(1-a^2)} aG_p \geq 2^{\frac12+\frac1p}   \sqrt{1- e^{-1/p} G_p^{-2}} e^{-1/2p}
\end{equation}
since the function $a\mapsto a\sqrt{1-a^2}$ is decreasing on $[1/\sqrt{2},1]$. Let $b=\frac1p \in[\frac{1}{26},\frac{1}{4}]$. The right-hand side of \eqref{eq:61} is at least 1 if and only if 
\begin{equation}
\log G_{1/b}^2 +b + \log\big(1-\tfrac12 \big(\tfrac{e}{4}\big)^b\big) \geq 0.
\end{equation}
Denote by $k(b)$ the resulting expression. The function $b\mapsto \log G_{1/b}^2$ is concave, as shown in the proof of Lemma \ref{lem:large-coordinate-1} and a derivative computation shows that the same holds for the function $b\mapsto b+\log(1-\tfrac12 (e/4)^b)$. Therefore, $k(b) \geq \min\{k(1/4),k(1/26)\} \approx 0.0001\ldots >0$.

If on the other hand $q_p^*>2$, i.e. $aG_p > e^{-1/2p}$, then
\begin{equation}
\sup_{q\geq2} \alpha_p(q) = \alpha_p(q_p^*) = \frac{\sqrt{2(1-a^2)} }{(ep)^{1/p} aG_p (\log(1/aG_p))^{1/p}} =\sqrt{ 2\Big(e^{\frac{2}{pq_p^\ast}} - G_p^{-2} \big) } \Big(\frac{q_p^\ast}{e}\Big)^{1/p}.
\end{equation}
This quantity is at least one if and only if the function
\begin{equation}
H(p,q) \eqdef \log G_p^2 +\log \bigg( e^{\frac{2}{pq}}-\frac12\Big(\frac{e}{q}\Big)^{2/p} \bigg)
\end{equation}
satisfies $H(p,q_p^\ast)\geq0$. Fixing $q>2$, one can again use the concavity of $b\mapsto \log G_{1/b}^2$ and a straightforward derivative computation to show that the function $b \mapsto H(1/b,q)$ is concave, therefore $H(p,q) \geq \min\{H(4,q), H(26,q)\}$ for $4\leq p\leq26$. Using $e^x\geq x+1$, we can bound
\begin{equation}
e^{\frac{2}{pq}} \geq e^{\frac1p}  \big(1+\tfrac{2}{pq}-\tfrac1p\big),
\end{equation}
therefore
\begin{equation} \label{eq:qmin0}
H(p,q) \geq \log G_p^2 +\log \bigg( e^{\frac1p}  \Big(1+\frac{2}{pq}-\frac1p\Big) -\frac12\Big(\frac{e}{q}\Big)^{2/p}\bigg).
\end{equation}
Since $p>2$, the function inside the logarithm is convex in the variable $1/q$ and attains its minimum at $q_{\mr{min}} = (2 e^{-1/p})^{\frac{p}{p-2}}>2$. If $p=4$,  \eqref{eq:qmin0} gives $H(4,q_{\mathrm{min}}) > 0.07$ whereas if $p=26$,  $H(26,q_{\mr{min}}) > 0.0001$. This completes the proof.
\end{proof}

\begin{proof} [Proof of Lemma \ref{lem:weighted-gaussian-comparison}]
For $\psi(s)=e^{-c_ps^2}$ there is nothing to prove, since the left-hand side coincides with the term on the right without the factor $\frac65$ and the equality follows by a change of variables. Consider $\psi=\gamma_p$.  By Proposition \ref{prop:Qp} or K\"onig's inequality \eqref{eq:Konig-envelope}, for every $s\in[0,3]$,
\begin{equation} \label{eq:small-s}
\big|\gamma_p(s)\big|=\gamma_p(s)\le e^{-c_ps^2},
\end{equation}
where the first equality holds since $\gamma_p(s)>0$ for $s\in[0,3]$ as was shown in the proof of Corollary \ref{cor:Qp}. On the other hand, for $s\geq3$, by \cite[Lemma~4.3]{Kon25},
\begin{equation} 
\big|\gamma_p(s)\big| \leq 0.2337 < \frac{4}{17}
\end{equation}
and combining this with \eqref{eq:bp-1-time}, we get that for $s\geq3$,
\begin{equation} \label{eq:large-s}
\big|\gamma_p(s)\big| \leq \min\Big\{\frac{4}{17}, \frac{1}{\Gamma(1+\frac1p) s}\Big\}.
\end{equation} 
Hence
\begin{equation}
\begin{split}
\int_0^\infty s^4|\gamma_p(s)|^u\, \diff s & \stackrel{\eqref{eq:small-s}\wedge\eqref{eq:large-s}}{\leq} \int_0^3 s^4 e^{-uc_ps^2}\,\diff s + \Big(\frac{4}{17}\Big)^{u-6} \frac{1}{\Gamma(1+\tfrac1p)^6}\int_3^\infty \, \frac{\diff s}{s^2}
\\ & \leq\int_0^\infty s^4 e^{-uc_ps^2}\,\diff s + \Big(\frac{4}{17}\Big)^{u-6} \frac{1}{3\Gamma(1+\tfrac1p)^6}.
\end{split}
\end{equation}
The ratio of the last tail bound to the desired Gaussian integral is equal to
\begin{equation}
u^{5/2}\Big(\frac{4}{17}\Big)^{u-6} \frac{8c_p^{5/2}}{9\sqrt{\pi}\Gamma(1+\tfrac1p)^6} \leq  8^{5/2}\Big(\frac{4}{17}\Big)^{2} \frac{8\cdot 6^{-5/2}}{9\sqrt{\pi}\Gamma(\frac{27}{26})^6} < \frac{1}{5},
\end{equation}
where we used that the expression on the left is decreasing in $u\geq8$, $c_p\leq \tfrac16$ (see, e.g., \cite[Lemma~2.4(b)]{Kon25}) and $\Gamma(1+\frac1p)\geq\Gamma(\frac{27}{26})$ for $p\geq 26$. This completes the proof of \eqref{eq:weighted-gaussian-comparison}.
\end{proof}

\begin{proof} [Proof of Lemma \ref{lem:3.3}]
We repeat the proof of Lemma \ref{lem:large-coordinate-2}. Let $M_p = \max\{T_p,G_p\}$. If $aM_p\geq1$, the function $\alpha_p(q) \eqdef \sqrt{2(1-a^2)}q^{1/p}(aM_p)^{q-1}$ has $\lim_{q\to\infty} \alpha_p(q)=\infty$. If $aM_p<1$, then the derivative $\alpha_p'$ vanishes at $q_p^* = \frac{1}{p\log(1/aM_p)}$. If $q_p^*\leq 2$, that is $aM_p \leq e^{-1/2p}$, then
\begin{equation} \label{eq:61-2}
\sup_{q\geq2} \alpha_p(q) = \alpha_p(2) = 2^{\frac12+\frac1p} \sqrt{(1-a^2)} aM_p \geq 2^{\frac12+\frac1p} \sqrt{1- e^{-1/p} M_p^{-2}} e^{-1/2p}
\end{equation}
since the function $a\mapsto a\sqrt{1-a^2}$ is decreasing on $[1/\sqrt{2},1]$. Let $b=\frac1p \in(0,\frac{1}{26}]$. Inequality \eqref{eq:lem-3.3} is satisfied provided that
\begin{equation}
\log M_{1/b}^2 +b + \log\big(1-\tfrac{1}{2e^{b^2/6}} \big(\tfrac{e}{4}\big)^b\big) \geq 0.
\end{equation}
Denote by $k(b)$ the resulting expression. Since $M_{1/b}\geq T_{1/b}=2^{\frac12-b}$, we have 
\begin{equation}
k(b)\geq \log\big(2(\tfrac{e}{4})^b-e^{-b^2/6}(\tfrac{e}{4})^{2b}\big),
\end{equation} which is positive if and only if $\big((\tfrac4e)^b-1\big)^2<1-e^{-b^2/6}$. Using $e^x-1\leq xe^x$ and $1-e^{-x}\geq x-\frac{x^2}{2}$, this follows from
\begin{equation}
\log\big(\tfrac4e\big)^2 \big(\tfrac4e\big)^{2b} < \tfrac16-\tfrac{b^2}{72};
\end{equation}
this holds for $b\leq\frac1{26}$, since the left-hand side is increasing and the right-hand side is decreasing in $b$, and at $b=\frac{1}{26}$ the left-hand side is roughly $0.154$ while the right-hand side is roughly $0.166$.

If on the other hand $q_p^*>2$, i.e. $aM_p > e^{-1/2p}$, then
\begin{equation}
\sup_{q\geq2} \alpha_p(q) = \alpha_p(q_p^*) = \frac{\sqrt{2(1-a^2)} }{(ep)^{1/p} aM_p (\log(1/aM_p))^{1/p}} =\sqrt{ 2\Big(e^{\frac{2}{pq_p^\ast}} - M_p^{-2} \big) } \Big(\frac{q_p^\ast}{e}\Big)^{1/p}.
\end{equation}
Inequality \eqref{eq:lem-3.3} is satisfied if and only if the function
\begin{equation}
H(p,q) \eqdef \log M_p^2 +\log \bigg( e^{\frac{2}{pq}}-\frac{1}{2e^{1/6p^2}}\Big(\frac{e}{q}\Big)^{2/p} \bigg)
\end{equation}
satisfies $H(p,q_p^\ast)\geq0$. Using $e^x\geq x+1$, we can bound
\begin{equation}
e^{\frac{2}{pq}} \geq e^{\frac1p} \big(1+\tfrac{2}{pq}-\tfrac1p\big),
\end{equation}
therefore, since $M_p\geq T_p$,
\begin{equation} \label{eq:qmin}
H(p,q) \geq \log T_p^2 +\log \bigg( e^{\frac1p} \Big(1+\frac{2}{pq}-\frac1p\Big) -\frac{1}{2e^{1/6p^2}}\Big(\frac{e}{q}\Big)^{2/p}\bigg).
\end{equation}
Since $p>2$, the function inside the logarithm is convex in the variable $1/q$ and attains its minimum at $q_{\mr{min}} = (2 e^{1/6p^2-1/p})^{\frac{p}{p-2}}$, where, with $b=\frac1p$,
$\frac{1}{q_{\mr{min}}}=\frac12\exp\big(-\frac{b\log(4/e)+b^2/6}{1-2b}\big)$. Plugging it in, the right-hand side of \eqref{eq:qmin} then becomes 
\begin{equation}
\log\big(2(\tfrac{e}{4})^b(1-b-\tfrac{1-2b}{q_{\mr{min}}})\big).
\end{equation} 
Using $e^{-x}\leq 1-x+\frac{x^2}{2}$ and $e^x\leq 1+x+\frac{x^2}{2}+\frac{x^3}{6}e^x$ for $x\geq0$, this is positive provided that
\begin{equation}
\frac{\big(\log\frac4e+\frac b6\big)^2}{2(1-2b)}+\tfrac{1}{2}\log^2\big(\tfrac4e\big)+\tfrac{b}{6}\log^3\big(\tfrac4e\big)\big(\tfrac4e\big)^b < \frac16 .
\end{equation}
This function is increasing in $b$ and at $b=\frac1{26}$ it is smaller than $0.159<\frac16$. 
\end{proof}


\medskip 

\noindent {\bf Disclosure of AI use.} During the preparation of this work, AI was exclusively (and extensively) used for numerical purposes, to test viable conjectures by plotting functions of a single real variable and for the numerical evaluation of series, integrals and special functions. The only exceptions are the one-dimensional inequalities \eqref{eq:k2k3} and \eqref{eq:P(3)}, whose proofs were generated by AI. The entirety of the mathematical reasoning is due to the authors, as is the exposition of the obtained results and the bibliographic research. After the paper was completed, AI was used to proofread the manuscript for typos and grammatical errors.

\bibliographystyle{plain}
\bibliography{max-sec-p>2}

\end{document}